\documentclass[12pt]{article}%
\usepackage{graphicx}
\usepackage{amsmath}%
\usepackage{amsfonts}%
\usepackage{amssymb}
\newtheorem{theorem}{Theorem}[section]

\newtheorem{case}{Case}

\newtheorem{criterion}[theorem]{Criterion}
\newtheorem{definition}[theorem]{Definition}
\newtheorem{example}[theorem]{Example}

\newtheorem{lemma}[theorem]{Lemma}

\newtheorem{proposition}[theorem]{Proposition}
\newtheorem{remark}[theorem]{Remark}

\newenvironment{proof}[1][Proof]{\textbf{#1.} }{\ \rule{0.5em}{0.5em}}

\begin{document}

\author{Peter Scott\thanks{Partially supported by NSF grant DMS 034681}\\Mathematics Department, University of Michigan\\Ann Arbor, Michigan 48109, USA.\\email:pscott@math.lsa.umich.edu
\and Gadde A. Swarup\\Mathematics Department, University of Melbourne\\Parkville, Victoria 3052, Australia.\\email:gadde@ms.unimelb.edu.au}
\title{Canonical splittings of groups and 3-manifolds}
\maketitle

\begin{abstract}
We introduce the notion of a `canonical' splitting over $\mathbb{Z}$ or
$\mathbb{Z}\times\mathbb{Z}$ for a finitely generated group $G$. We show that
when $G$ happens to be the fundamental group of an orientable Haken manifold
$M$ with incompressible boundary, then the decomposition of the group
naturally obtained from canonical splittings is closely related to the one
given by the standard JSJ-decomposition of $M$. This leads to a new proof of
Johannson's Deformation Theorem.

\end{abstract}
\date{}

An important chapter in the structure theory of Haken $3$-manifolds is the
theory of the characteristic submanifold initiated by Waldhausen in
\cite{Waldhausen1}. This was accomplished in 1974-75 by Jaco and Shalen
\cite{JS} and Johannson \cite{J}. They defined the characteristic submanifold
$V(M)$ of a Haken $3$-manifold $M$ with incompressible boundary and proved it
to be unique up to isotopy. It has the properties that each component is
either a Seifert fibre space or is an $I$-bundle over a surface, and every
essential map of an annulus or torus into $M$ is properly homotopic into
$V(M)$. Johannson called this last statement the Enclosing Property. This
decomposition of $M$ into submanifolds has come to be called the
JSJ-decomposition. Johannson emphasised the importance of the
JSJ-decomposition by using it to prove his Deformation Theorem \cite{J}, which
gives a good description of the connection between Haken $3$-manifolds which
have isomorphic fundamental groups. More recently Sela \cite{Sela} (see also
Bowditch \cite{B}) showed that analogous algebraic decompositions exist for
word hyperbolic groups, and Rips and Sela extended aspects of this theory to
torsion free finitely generated groups \cite{RipsSela}. Dunwoody and Sageev
\cite{DS} showed that such algebraic decompositions exist for all finitely
presented groups, and Fujiwara-Papasoglu \cite{FP} further extended this
theory. However none of these group theoretic decompositions yields the usual
JSJ-decomposition when restricted to fundamental groups of orientable Haken
manifolds with incompressible boundary.

In this paper, we describe the standard JSJ-decomposition of an orientable
Haken $3$-manifold in terms that can be translated into group theory and we
then give a purely algebraic description of it. This leads to a natural
algebraic proof of Johannson's Deformation Theorem \cite{J}. The papers
\cite{NS}, \cite{LS}, and \cite{SS} represent our earlier work in this
programme. Our description of the standard JSJ-decomposition is the following.
Instead of trying to characterize the characteristic submanifold, we
concentrate on its frontier which consists of essential annuli and tori. Let
$M$ be an orientable Haken manifold with incompressible boundary. We call a
two-sided embedded essential annulus or torus $S$ in $M$ \textit{canonical} if
any two-sided immersed essential annulus or torus $T$ in $M$ can be homotoped
away from $S$. The Enclosing Property of the characteristic submanifold $V(M)$
of $M$ implies that any component of the frontier of $V(M)$ is canonical. We
will see in section 1 that, up to isotopy, these are the only canonical
surfaces in $M$. Without assuming this, we can consider the collection of all
isotopy classes of canonical surfaces in $M$ and pick one surface from each
isotopy class. Up to isotopy, any pair of these surfaces will be disjoint. By
choosing a Riemannian metric on $M$ and choosing each surface to be least
area, it follows that we can choose all these surfaces to be disjoint. Thus we
obtain a family of disjoint embedded essential annuli and tori in $M$ such
that no two of them are parallel. There is an upper bound to the number of
such surfaces in $M$, which is called the Haken number of $M$. In particular,
it follows that this must be a finite family, and we call it a JSJ-system for
$M$. Such a system is unique up to isotopy. We will see in section 1 that a
JSJ-system in $M$ can be obtained from the frontier of the characteristic
submanifold $V(M)$ of $M$ by choosing one surface from each maximal collection
of parallel surfaces. The reader is warned that in \cite{NS}, Neumann and
Swarup defined the term canonical in a slightly different way. They called $S$
canonical if any \emph{embedded} essential annulus or torus in $M$ can be
isotoped away from $S$. The two concepts are not equivalent. For example, let
$M$ be constructed by gluing two solid tori $V$ and $W$ along an annulus $A$
which is embedded in their boundaries with degree greater than $1$ in each.
Thus $M$ is a Seifert fibre space and $A$ is an essential annulus in $M$. Now
$M$ contains no essential embedded tori and the only embedded essential
annulus in $M$ is $A$. Thus $A$ is obviously canonical in the sense of Neumann
and Swarup. But $M$ does contain singular essential tori, and any such must
intersect $A$, so that $A$ is not canonical in our sense.

In order to describe an algebraic analogue of a JSJ-system, we need to have
algebraic versions of the ideas of embedded surfaces, of immersed surfaces and
of the property of being canonical. The analogue of a two-sided incompressible
surface embedded in a $3$-manifold is a splitting of a group along a surface
group, which in our case will be a splitting of $\pi_{1}(M)$ along
$\mathbb{Z}$ or $\mathbb{Z}\times\mathbb{Z}$. The notion of a $\pi_{1}%
$-injective immersed surface $F$ in a $3$-manifold $M$ corresponds to a
$H$-almost invariant subset of $G$, where $G$ denotes $\pi_{1}(M)$ and $H$
denotes the image of $\pi_{1}(F)$ in $G$. Such sets have been studied by Scott
in \cite{Scott:TorusTheorem} and \cite{Scott:Strong}, Dunwoody and others (see
the references in \cite{SS}). The notion of homotopic disjointness of immersed
surfaces can be formulated in terms of the intersection numbers developed by
Freedman, Hass and Scott in \cite{FHS}. It follows from that paper that two
immersed surfaces can be homotoped to have disjoint images if and only if
their intersection number equals zero. In \cite{Scott:Intersectionnumbers},
Scott showed that this notion of intersection number can be naturally extended
to one for almost invariant sets. These ideas were further extended in
\cite{SS} where we showed that almost invariant sets with self-intersection
number zero naturally produce splittings, and that splittings with
intersection number zero can be made ``disjoint''. (In \cite{SS}, the idea of
disjointness of splittings is formulated precisely in terms of graphs of
groups. See the end of section 1 of this paper.) Thus, we have all the notions
needed to reformulate the above description of a JSJ-system in a $3$-manifold
$M$ in a purely algebraic fashion.

One of the important aspects of the JSJ-decomposition for $3$-manifolds is the
description of the pieces obtained by splitting $M$ along the JSJ-system.
These are either simple or fibred or both and the characteristic submanifold
$V(M)$ is essentially the union of the fibred pieces (see section
\ref{deformationtheoremchapter}). Since the Enclosing Property implies that
any immersed essential annulus or torus in $M$ can be properly homotoped into
$V(M)$, this may be taken as a description of all immersed essential annuli
and tori in $M$, in particular of all embedded essential annuli and tori. It
is these aspects that the group theoretic generalizations (except the one for
word hyperbolic groups) do not capture fully. Apart from $I$-bundle pieces,
the characteristic submanifold consists of Seifert fibred pieces. The role
played by those Seifert pieces which meet the boundary of $M$ (we will call
them \textit{peripheral}) is rather different from that played by the other
(interior) Seifert pieces. In all the group theoretic JSJ-decompositions, the
peripheral Seifert fibred pieces are further decomposed, and the interior ones
may or may not appear depending on which theory is being considered; but all
the theories capture the analogues of the $I$-bundle pieces. The general
difficulty seems to come from the fact that the both annuli and tori are
needed to describe the standard JSJ-decomposition whereas the group theoretic
decompositions so far work well only with $\mathbb{Z}$ or $\mathbb{Z}%
\times\mathbb{Z}$ splittings, but not both, and one type of intersection which
we call strong crossing in \cite{SS}. Perhaps it may not be possible to fully
capture all aspects of the standard JSJ-decomposition in the case of groups.

In section 1, we discuss background material on $3$-manifolds and group
theory. In section \ref{canonicalsplittings}, we study the analogues of
canonical annuli and tori, namely canonical $\mathbb{Z}$ or $\mathbb{Z}%
\times\mathbb{Z}$ splittings of $G=\pi_{1}(M)$. If there is a true analogue
for groups of the standard JSJ-decomposition for manifolds, one might hope
that canonical annuli and tori correspond to canonical splittings. This turns
out to be almost true. The exceptions are the canonical tori which have
peripheral Seifert fibred spaces on both sides. However, we are also able to
give a completely algebraic characterisation of these canonical tori. Since
this characterisation is algebraic, it follows that the JSJ-systems of any two
Haken manifolds with isomorphic fundamental groups correspond. In section
\ref{deformationtheoremchapter}, we use this fact to give a `natural'
algebraic proof of Johannson's Deformation Theorem \cite{J}. This result
states that if $M$ and $N$ are Haken $3$-manifolds with incompressible
boundary, and $f:M\rightarrow N$ is a homotopy equivalence, then $f$ is
homotopic to a map $g$ such that $g(V(M))\subset V(N)$, the restriction of $g$
to $V(M)$ is a homotopy equivalence from $V(M)$ to $V(N)$, and the restriction
of $g$ to $\overline{M-V(M)}$ is a homeomorphism onto $\overline{N-V(N)}$.
Most of this result follows immediately from our algebraic description of the
JSJ-decomposition. The proof of the last part of the Deformation Theorem which
asserts the existence of a homeomorphism on the complement of $V(M)$ is
carried out in section \ref{deformationtheoremchapter}. This follows the ideas
of Swarup \cite{Swarup:Deformation} and of Jaco \cite{JacoCBMSNotes}, and
corrects an omission in \cite{JacoCBMSNotes}. We remark that the proof in
\cite{Swarup:Deformation} is `algebraic' and is one of the reasons for our
algebraic point of view. The fact that topologically canonical tori need not
correspond to algebraically canonical splittings indicates one of the
difficulties in extending fully the standard JSJ-decomposition to groups.
Since these exceptions have also an algebraic description, there seems to be
some hope of carrying out some JSJ-decomposition for special classes of
groups. The ideas here can also be used to obtain the decompositions of Sela
\cite{Sela} and Bowditch \cite{B} for word hyperbolic groups.

\section{Background material in topology and algebra}

Recall that we say that a two-sided embedded essential annulus or torus $S$ in
$M$ is canonical if any two-sided immersed essential annulus or torus $T$ in
$M$ can be homotoped away from $S$. It will be convenient to reformulate the
definition of canonical in the following way, using the definition of
intersection number of surfaces in a $3$-manifold which was given in
\cite{FHS}.

\begin{definition}
An essential annulus or torus embedded in $M$ is \textsl{canonical} if it has
intersection number zero with any (possibly singular) essential annulus or
torus in $M$.
\end{definition}

Recall that a JSJ-system in $M$ is a family of disjoint embedded annuli and
tori ${S}_{1},...,S_{n}$ in $M$ which contains one representative of each
isotopy class of canonical annuli and tori. Such a system exists and is unique
up to isotopy. We now describe the relationship between the JSJ-system in $M$
and the characteristic submanifold $V(M)$ of $M$. Let $\mathcal{F}$ denote the
family of annuli and tori in $M$ obtained from the frontier $frV(M)$ of $V(M)$
by choosing one surface from each maximal collection of parallel surfaces.

\begin{proposition}
\label{JSJ-systemequalsF}Let $M$ be an orientable Haken $3$-manifold with
incompressible boundary. Then $\mathcal{F}$ is a JSJ-system of $M$.
\end{proposition}

\begin{proof}
The Enclosing Property of $V(M)$ shows that any component of $frV(M)$ is
canonical. Thus $\mathcal{F}$ consists of non-parallel canonical surfaces and
so is contained in a JSJ-system. Now any embedded essential annulus or torus
can be isotoped to lie in $V(M)$. Thus the proof of the proposition reduces to
showing that if a component $X$ of $V(M)$ contains an embedded annulus or
torus $S$ which is essential in $M$, and is not parallel to a component of
$frV(M)$, then $S$ is not canonical. This means that there is some essential
annulus or torus $S^{\prime}$ in $M$ which has non-zero intersection number
with $S$, or equivalently that $S^{\prime}$ cannot be properly homotoped apart
from $S$. This can be proved by careful consideration of the possibilities
using the fact that $X$ is a Seifert fibre space or an $I$-bundle and that $S$
can be isotoped to be vertical in $X$, i.e. to be a union of fibres in $X$.
However, this consideration has essentially been done already in \cite{NS}, so
we will use their work to avoid a case by case discussion.

In \cite{NS}, Neumann and Swarup defined what they called the $W$-system of
$M$. This consisted of a maximal system of disjoint annuli and tori ${S}%
_{1},...,S_{n}$ in $M$ so that no two of the $S_{i}$ are parallel, and each
$S_{i}$ is canonical in their sense, which means that any \emph{embedded}
essential annulus or torus in $M$ can be isotoped away from each $S_{i}$. We
would like to remark that the $W$-system is another natural concept that can
be generalized to groups. The $W$-system is unique up to isotopy and they
showed that the $W$-system contains $\mathcal{F}$. By the definition of the
$W$-system, any essential annulus or torus disjoint from and not parallel into
the $W$-system intersects some other embedded annulus or torus in an essential
way and so does not belong to a JSJ-system. The only surfaces which are in the
$W$-system but not in $\mathcal{F}$ are the special annuli described in Lemma
2.9 and Figure 5 of \cite{NS}. It is easy to see that in these cases (which
essentially arise from Seifert bundles over twice punctured discs), there is
an immersed annulus intersecting a special annulus in an essential way. It
follows that these special annuli also do not belong to the JSJ-system, so
that the JSJ-system consists precisely of the $W$-system with these special
annuli removed, which is the same as $\mathcal{F}$.
\end{proof}

In this paper, we are restricting to the orientable case for simplicity, but
it is easy to extend our analysis to the general case as follows. Let $M$ be a
non-orientable Haken $3$-manifold, let $M_{1}$ denote its orientable double
cover and consider the canonical family $\mathcal{F}_{1}$ of annuli and tori
in $M_{1}$. In Theorem 8.6 of \cite{Meeks-Scott}, it is shown that the tori in
$\mathcal{F}_{1}$ can be chosen to be invariant under the covering involution
unless $M_{1}$ is a bundle over the circle with fibre the torus and a
hyperbolic gluing map. (Note that the full power of the results of
\cite{Meeks-Scott} is not needed here. One need only know that any involution
of $T\times I$ respects some product structure, and this was first proved by
Kim and Tollefson \cite{KimTollefson}.) Although it was not discussed in
\cite{Meeks-Scott}, the same argument shows that $\mathcal{F}_{1}$ itself can
be chosen to be invariant. The only possible new exceptional case would be if
$M_{1}$ is a bundle over the circle with fibre the annulus. But any such
manifold is a Seifert fibre space, so that $\mathcal{F}_{1}$ is empty and is
trivially invariant in this case. In particular, if $M$ has non-empty
boundary, we can always choose $\mathcal{F}_{1}$ to be invariant under the
covering involution. The image of $\mathcal{F}_{1}$ in $M$ is a collection of
embedded essential annuli, tori, Moebius bands and Klein bottles some of which
may be one-sided. If one of these surfaces is one-sided, we replace it by the
boundary of a thin regular neighbourhood. The resulting collection of
two-sided surfaces in $M$ is denoted by $\mathcal{F}$. As in the orientable
case, the surfaces in $\mathcal{F}$ are characterised by the property that
they are two-sided and essential, and have intersection number zero with every
essential map of an annulus or torus into $M$.

We will need the following facts about the canonical decomposition of $M$
obtained by cutting along the JSJ-system $\mathcal{F}$. The manifolds so
obtained are called the canonical pieces of the decomposition. If a Seifert
fibred canonical piece $\Sigma$ contains canonical annuli in its frontier,
then $\Sigma$ admits a Seifert fibration for which all these canonical annuli
are vertical in $\Sigma$. Note that if $\Sigma$ is not closed, it has a unique
Seifert fibration unless it is a $I$-bundle over the Klein bottle. We will
need to give special consideration to any canonical torus $T$ of $M$ such that
each of the canonical pieces of $M$ adjacent to $T$ is a Seifert fibre space
which is peripheral, i.e. meets the boundary of $M$. We will say that such a
canonical torus is \textit{special}. Note that it is possible that a single
Seifert fibred canonical piece of $M$ meets both sides of a special canonical torus.

In \cite{Scott:Intersectionnumbers}, Scott defined intersection numbers for
splittings of groups and more generally for almost invariant subsets of coset
spaces of a group. See \cite{Scott:Intersectionnumbers} and \cite{SS} for a
detailed discussion of the connection of the algebraic idea of intersection
number with the topological idea. In particular, let $G$ denote $\pi_{1}(M)$,
pick a basepoint in $M$ and pick finitely many simple loops representing
generators of $G$, which are disjoint except at the basepoint. Let $\Gamma$
denote the pre-image of the union of these loops in the universal cover
$\widetilde{M}$ of $M$. Once we have chosen a vertex of this graph as our
basepoint in $\widetilde{M}$, we can identify the vertices of $\Gamma$ with
$G$, and can identify $\Gamma$ itself with the Cayley graph of $G$ with
respect to the given generating set. If $H$ is a subgroup of $G$, the quotient
graph $H\backslash\Gamma$ is embedded in the quotient manifold $H\backslash
\widetilde{M}$, which we will usually denote by $M_{H}$, and the vertices of
$H\backslash\Gamma$ are naturally identified with the cosets $Hg$ of $H$ in
$G$. Let $F$ denote an annulus or torus and let $f:F\rightarrow M$ denote an
essential map. This means that $f$ is proper and $\pi_{1}$-injective, and that
$f$ is not properly homotopic into the boundary $\partial M$ of $M$. Let $H$
denote the subgroup $f_{\ast}\pi_{1}(F)$ of $G$. Let $M_{F}$ denote the cover
of $M$ such that $\pi_{1}(M_{F})=H$, and consider the lift of $f$ into $M_{F}%
$. For convenience, suppose that this lift is an embedding. (This is true up
to homotopy and is automatically true if $f$ is least area by \cite{FHS}.) Now
let $P$ and $P^{\ast}$ denote the two submanifolds of $M_{F}$ into which $F$
cuts $M_{F}$. Thus the frontier in $M_{F}$ of each of $P$ and $P^{\ast}$ is
exactly $F$. As $f$ is essential, neither $P$ nor $P^{\ast}$ is compact. For
if $P$ were compact, it would have to be homeomorphic to $F\times I$ and so
the lift of $F$ would be properly homotopic into $\partial M_{F}$,
contradicting our assumption that $f$ is not properly homotopic into $\partial
M$. We associate to $P$ the subset of $H\backslash G$ consisting of those
vertices of $H\backslash\Gamma$ which lie in the interior of $P$. Because the
frontier of $P$ is compact, this subset has finite coboundary in
$H\backslash\Gamma$, and so is an almost invariant subset of $H\backslash G$.
As neither $P$ nor $P^{\ast}$ is compact, this associated almost invariant set
is non-trivial, which means that it is infinite with infinite complement.
(Note that we are interested only in the almost equality class of this almost
invariant subset, so it does not matter if some of the vertices of
$H\backslash\Gamma$ lie on $F$.) For our purposes in this paper, we want to
think of $P$ as being much the same as its associated almost invariant set.
However, in order to avoid confusion, we will not refer to $P$ as an almost
invariant set. Instead we make the following definitions.

\begin{definition}
A submanifold $P$ (which need not be connected) of a $3$-manifold $W$ is an
\textsl{ending submanifold} of $W$ if $P$ and the closure of $W-P$ are not
compact, but the frontier of $P$ is compact. Two ending submanifolds $P$ and
$Q$ of $W$ are \textsl{almost equal} if their symmetric difference $P-Q\cup
Q-P$ is bounded.
\end{definition}

Using the natural map which we have just defined from ending submanifolds of
$M_{H}$ to non-trivial almost invariant subsets of $H\backslash G$, leads to
the following definition of the intersection number of two ending submanifolds.

\begin{definition}
Let $M$ be a compact $3$-manifold with fundamental group $G$, let $H$ and $K$
be subgroups of $G$, and let $P\subset M_{H}$ and $Q\subset M_{K}$ be ending
submanifolds. Let $X$ and $Y$ denote the pre-images of $P$ and $Q$
respectively in the universal cover $\widetilde{M}$ of $M$. Then
$Y\;$\textsl{crosses} $X$ if and only if each of the four sets $X\cap Y$,
$X\cap Y^{\ast}$, $X^{\ast}\cap Y$, $X^{\ast}\cap Y^{\ast}$ has unbounded
image in $M_{H}$. The \textsl{intersection number} of $P$ and $Q$ is the
number of double cosets $HgK$ such that $gY$ crosses $X$.
\end{definition}

It follows from \cite{Scott:Intersectionnumbers} that crossing is symmetric,
so that $Y$ crosses $X$ if and only if $X$ crosses $Y$. Also the intersection
number of $P$ and $Q$ is unaltered if either set is replaced by its complement
or by an almost equal ending submanifold.

Our previous discussion shows how to associate an ending submanifold $P$ of
$M_{F}$ to an essential annulus or torus in $M$. This submanifold is well
defined up to complementation (we could equally well choose $P^{\ast}$ in
place of $P$), and up to almost equality. Note that $P$ can be defined
perfectly well even if the lift of $f$ is not embedded. For the complement of
the image of the lift will have two unbounded components, and we choose one of
them as $P$. We will need the following two definitions.

\begin{definition}
Given a splitting of $G$ over $H$, where $H$ is isomorphic to $\mathbb{Z}$ or
$\mathbb{Z}\times\mathbb{Z}$, we will define this splitting to be
\textsl{algebraically canonical} if the corresponding almost invariant subset
of $H\backslash G$ has zero intersection number with any almost invariant
subset of any $K\backslash G$, for which $K$ is isomorphic to $\mathbb{Z}$ or
$\mathbb{Z}\times\mathbb{Z}$.
\end{definition}

\begin{definition}
Given a splitting of $G$ over $H$, where $H$ is isomorphic to $\mathbb{Z}$ or
$\mathbb{Z}\times\mathbb{Z}$, we will define this splitting to be $\mathbb{Z}%
$-\textsl{canonical} if the corresponding almost invariant subset of
$H\backslash G$ has zero intersection number with any almost invariant subset
of any $K\backslash G$, for which $K$ is isomorphic to $\mathbb{Z}$.
\end{definition}

We have now defined what it means for an essential annulus or torus in $M$ to
be topologically canonical and for a splitting of $G$ to be algebraically
canonical or $\mathbb{Z}$-canonical. In order to compare these algebraic and
topological ideas, we reformulate them all in terms of ending submanifolds.
Note that a splitting of $G=\pi_{1}(M)$ over $H$ has an associated almost
invariant subset of $H\backslash G$, and we can find an ending submanifold $P$
of $M_{H}$ whose associated almost invariant set is exactly this. We will say
that $P$ is associated to the splitting. Of course, $P$ is not unique but any
two choices for $P$ will be almost equal or almost complementary. We have the
following statements.

\begin{criterion}
\label{criterion}

\begin{enumerate}
\item A splitting of $G$ over $H$, for which $H$ is isomorphic to $\mathbb{Z}$
or to $\mathbb{Z}\times\mathbb{Z}$, is algebraically canonical if and only if
the associated ending submanifold $P$ of $M_{H}$ has intersection number zero
with every ending submanifold $Q$ of every $M_{K}$, for which $K\;$is
isomorphic to $\mathbb{Z}$ or to $\mathbb{Z}\times\mathbb{Z}$.

\item A splitting of $G$ over $H$, for which $H$ is isomorphic to $\mathbb{Z}$
or to $\mathbb{Z}\times\mathbb{Z}$, is $\mathbb{Z}$-canonical if and only if
the associated ending submanifold $P$ of $M_{H}$ has intersection number zero
with every ending submanifold $Q$ of every $M_{K}$, for which $K\;$is
isomorphic to $\mathbb{Z}$.

\item An essential annulus or torus $F$ embedded in $M$, and carrying the
group $H$, is topologically canonical if and only if the associated ending
submanifold $P$ of $M_{H}$ has intersection number zero with every ending
submanifold $Q$ of every $M_{K}$, such that $K\;$is isomorphic to $\mathbb{Z}$
or to $\mathbb{Z}\times\mathbb{Z}$ and $Q$ is associated to an essential
annulus or torus in $M$.
\end{enumerate}
\end{criterion}

Note that there is a natural topological idea corresponding to the algebraic
idea of a splitting being $\mathbb{Z}$-canonical. This is that an embedded
essential torus or annulus in $M$ is \textit{annulus-canonical} if it has zero
intersection number with every essential annulus in $M$. However, we will not
need to use this idea in this paper.

Recall from the introduction that if we have a collection of canonical
embedded annuli and tori in $M$, they can be isotoped to be disjoint. We
discussed the algebraic analogue in \cite{SS}, and will briefly describe the
results proved there. Let $S_{1},\ldots,S_{n}$ be a family of disjoint
essential embedded annuli or tori in a $3$-manifold $M$ with fundamental group
$G$, so that each $S_{i}$ determines a splitting of $G$. Together they
determine a graph of groups structure on $G$ with $n$ edges. A collection of
$n$ splittings of a group $G$ is said to be \textit{compatible} if $G$ can be
expressed as the fundamental group of a graph of groups with $n$ edges, such
that, for each $i$, collapsing all edges but the $i$-th yields the $i$-th
splitting of $G.$ The splittings are \textit{compatible up to conjugacy} if
collapsing all edges but the $i$-th yields a splitting of $G$ which is
conjugate to the $i$-th given splitting. Clearly disjoint essential annuli and
tori in $M$ define splittings of $G$ which are compatible up to conjugacy. The
following result is Theorem 2.5 of \cite{SS}.

\begin{theorem}
\label{splittingsarecompatibleiffzerointersectionnumber}\textit{Let }%
$G$\textit{\ be a finitely generated group with }$n$\textit{\ splittings over
finitely generated subgroups. This collection of splittings is compatible up
to conjugacy if and only if each pair of splittings has intersection number
zero. Further, in this situation, the graph of groups structure on }$G$ has a
unique underlying graph, and the edge and vertex groups are unique up to conjugacy.
\end{theorem}

Now let $G$ denote the fundamental group of a Haken $3$-manifold, and consider
a finite collection of algebraically canonical splittings of $G$ over
subgroups which are isomorphic to $\mathbb{Z}$ or to $\mathbb{Z}%
\times\mathbb{Z}$. Certainly any two such splittings have zero intersection
number, so together they yield a graph of groups structure for $G$ such that,
for each $i$, collapsing all edges but the $i$-th yields a splitting of $G$
which is conjugate to the $i$-th given splitting$.$ Further this graph is unique.

\section{The main results}

\label{canonicalsplittings}

The aim of this section is to relate the topological and algebraic ideas of
being canonical which we have introduced in the previous section.

Our major result is the following.

\begin{theorem}
\label{nonspecialsurfacesandcanonicalsplittings}Let $M$ be an orientable Haken
$3$-manifold with incompressible boundary, and let $\Phi$ denote the natural
map from the set of isotopy classes of embedded essential annuli and tori in
$M$ to the set of conjugacy classes of splittings of $G=\pi_{1}(M)$ over a
subgroup isomorphic to $\mathbb{Z}$ or to $\mathbb{Z}\times\mathbb{Z}$.

\begin{enumerate}
\item If $\sigma$ is an algebraically canonical splitting of $G$ over
$\mathbb{Z}$ or $\mathbb{Z}\times\mathbb{Z}$, then there is a canonical
annulus or torus $F$ in $M$, such that $\Phi(F)=\sigma$.

\item If $F$ is a canonical annulus or torus in $M$, then either $\Phi(F)$ is
an algebraically canonical splitting of $G$ over $\mathbb{Z}$ or
$\mathbb{Z}\times\mathbb{Z}$, or $F$ is a special canonical torus in $M$.
\end{enumerate}
\end{theorem}

\begin{remark}
We will show in Lemma \ref{specialtorusdoesnotimplycanonicalsplitting} that if
$F$ is a special canonical torus in $M$, then $\Phi(F)$ is never algebraically
canonical. Thus $\Phi$ induces a bijection between the collection of isotopy
classes of non-special canonical annuli and tori in $M$ and the collection of
conjugacy classes of algebraically canonical splittings of $G$. (Recall that,
by definition, the word special only applies to tori.)
\end{remark}

We prove part 1) of Theorem \ref{nonspecialsurfacesandcanonicalsplittings} in
Lemmas \ref{canonical Zsplittingimpliescanonicalannulus} and
\ref{canonicalZxZsplittingimpliescanonicaltorus}, and prove part 2) in Lemmas
\ref{canonicalannulusimpliescanonicalsplitting} and
\ref{nonspecialtorusimpliescanonicalsplitting}.

Note that if a canonical splitting of $G$ is induced by an embedded annulus or
torus $F$, it is obvious that $F$ is topologically canonical, by Criterion
\ref{criterion}. The difficulty is that many splittings of $G$ over
$\mathbb{Z}$ or $\mathbb{Z}\times\mathbb{Z}$ are not induced by a splitting
over a connected surface. For a very simple example of this phenomenon in one
dimension less, let $\Sigma$ denote the disc with two holes, i.e. the compact
planar surface with three boundary components. There are only six essential
arcs in $\Sigma$ each of which defines a splitting of $G=\pi_{1}(\Sigma)$ over
the trivial group. (Three of these arcs join distinct components of
$\partial\Sigma$, and the other three do not.) But there are infinitely many
non-conjugate splittings of the free group of rank two over the trivial group,
corresponding to the infinitely many choices of generators. Thus all of these
splittings except for six are not induced by a connected $1$-manifold in
$\Sigma$. For a $3$-dimensional example, we simply take the product of
$\Sigma$ with the circle $S^{1}$. This will yield a compact $3$-manifold
$\Sigma\times S^{1}$ for which infinitely many splittings over $\mathbb{Z}$
are not represented by a connected surface.

To start the proof, we consider general ending submanifolds in a cover $M_{H}$
of $M$, with $H$ infinite cyclic. We know that the boundary of $M_{H}$ is
incompressible, and so consists of planes and open annuli. If $A_{1}$ and
$A_{2}$ are annulus components of $\partial M_{H}$, there is an essential loop
on $A_{1}$ homotopic to a loop on $A_{2}$. The Annulus Theorem
\cite{Scott:TorusTheorem} then implies that there is an embedded
incompressible annulus in $M_{H}$ joining $A_{1}$ and $A_{2}$. It follows that
all the annulus components of $\partial M_{H}$ carry the same subgroup of $H$,
and hence so also does any properly embedded compact incompressible annulus in
$M_{H}$. It can be shown that $M_{H}$ compactifies to a solid torus, but we
will not need to use this fact.

\begin{lemma}
\label{Zgeneralendingsubmanifold}Any ending submanifold $P$ of $M_{H}$, for
which $H$ is isomorphic to $\mathbb{Z}$, is almost equal to one with the
property that the frontier $dP$ of $P$ consists only of embedded essential
annuli, and distinct boundary circles of $dP$ lie in distinct boundary
components of $M_{H}$.
\end{lemma}

\begin{proof}
Perform surgery on the frontier $dP$ of $P$ to make it $\pi_{1}$-injective in
$M_{H}$. This adds to $P$ or subtracts from $P$ compact sets, namely 2-handles
attached to $dP$, so that, by changing $P$ by a compact set, we can arrange
that $P$ has $\pi_{1}$-injective frontier $dP$. As $\pi_{1}(M_{H})=H$ is
isomorphic to $\mathbb{Z}$, the components of $dP$ can only be annuli, discs
or spheres. As $M_{H}$ is irreducible, any sphere bounds a ball, so that any
spheres in $dP$ can be removed by adding this ball to $P$ or removing it, as
appropriate. As $M_{H}$ has incompressible boundary, if $D$ is a disc properly
embedded in $M_{H}$ with boundary $C$, then $C$ bounds a disc $D^{\prime}$ in
$\partial M_{H}$. Thus $D\cup D^{\prime}$ is a sphere in $M_{H}$ and so bounds
a ball. As before, any discs in $dP$ can be removed by adding this ball to $P$
or removing it, as appropriate. Similarly inessential annuli can be removed
from $dP$. Thus we eventually arrive at a stage where $dP$ consists of a
finite number of compact essential annuli in $M_{H}$. If two boundary circles
of $dP$ lie in a single component $S$ of $\partial M_{K}$, they will cobound a
compact annulus $A$ in $S$, and we can alter $P$ by adding or removing a
regular neighbourhood of $A$. This reduces the number of boundary circles of
$dP$, so by repeating, we can ensure that distinct boundary circles of $dP$
lie in distinct boundary components of $M_{K}$ as required. (This move may
produce an inessential annulus component of $dP$, in which case we can remove
it as above.) This completes the proof of the lemma.
\end{proof}

The following technical result will be extremely useful in what follows.

\begin{lemma}
\label{annuluscutsintotwoinfinitepieces}Let $A_{1},\ldots,A_{k}$ be disjoint
compact incompressible annuli in $M_{H}$ such that distinct boundary circles
of the $A_{i}$'s lie in distinct boundary components of $M_{H}$. Let $R$
denote the closure of a component of $M_{H}-\cup_{i=1}^{k}A_{i}$, and let $A$
denote an embedded incompressible annulus in $R$ which joins $A_{1}$ to
$A_{2}$. Then $A$ cuts $R$ into two unbounded components.
\end{lemma}

\begin{proof}
Recall that $H$ is infinite cyclic. This implies that $A$ carries a subgroup
of finite index in $\pi_{1}(R)$, so that $A$ must separate $R$ into two pieces
$V$ and $W$. Now let $S_{1},\ldots,S_{4}$ denote the four components of
$\partial M_{H}$ which contain the four boundary components of $A_{1}$ and
$A_{2}$. Thus each $S_{i}$ is an open annulus and $A_{1}$ and $A_{2}$ cut each
$S_{i}$ into two unbounded annuli. As each $S_{i}$ contains no other boundary
component of any $A_{j}$, this means that each of $V$ and $W$ contains two of
these unbounded annuli, and so must be unbounded. This completes the proof of
the lemma.
\end{proof}

\begin{lemma}
\label{Zendingsubmanifoldforcanonicalsplitting}Suppose that $P$ is an ending
submanifold of $M_{H}$, for which $H$ is isomorphic to $\mathbb{Z}$, such that
$P\;$is associated to a canonical splitting of $G$ over $H$. Then $P$ is
almost equal to an ending submanifold whose frontier consists of a single
essential annulus which carries $H$.
\end{lemma}

\begin{proof}
Lemma \ref{Zgeneralendingsubmanifold} tells us that we can arrange that $dP$
consists only of embedded essential annuli, and that distinct boundary circles
of $dP$ lie in distinct boundary components of $M_{H}$.

First we will show that $dP$ is connected. Otherwise, pick two components of
$dP$. There is a compact incompressible annulus $A$ embedded in $M_{H}$ and
joining these two components. It is possible that $A$ meets other components
of $dP$ in its interior, but we can isotop $A$ so that it meets each component
of $dP$ in essential circles only. Thus $A$ has a sub-annulus which joins
distinct components $A_{1}$ and $A_{2}$ of $dP$ and whose interior does not
meet $dP$. We will replace $A$ by this sub-annulus and will continue to call
this annulus $A$. Without loss of generality we can suppose that $A$ lies in a
component $Q$ of $P^{\ast}$, and cuts $Q$ into pieces $V$ and $W$. Lemma
\ref{annuluscutsintotwoinfinitepieces} shows that $V$ and $W$ must each be
unbounded. Let $P_{1}$ and $P_{2}$ denote the components of $P$ which contain
$A_{1}$ and $A_{2}$ respectively and let $R$ denote the submanifold $V\cup
P_{1}$ of $M_{H}$. Then all four of the intersection sets of $R$ and $R^{\ast
}$ with $P$ and $P^{\ast}$ are unbounded, so we have found an ending
submanifold of $M_{H}$ with non-zero intersection number with the ending
submanifold $P$, contradicting our hypothesis that $P\;$is associated to a
canonical splitting of $G$ over $H$. This contradiction shows that $dP$ must
be connected.

Next we show that $dP$ must carry $H$. Otherwise, $dP$ carries a proper
subgroup $K$ of finite index in $H$. Let $Q$ denote the pre-image of $P$ in
$M_{K}$, so that $dQ$ is the pre-image of $dP$ in $M_{K}$ and consists of at
least two embedded essential annuli. Note that distinct boundary circles of
$dQ$ must lie in distinct boundary components of $M_{K}$. Now we can apply the
argument of the preceding paragraph with $M_{K}$ in place of $M_{H}$ and $Q$
in place of $P$. This will yield an ending submanifold of $M_{K}$ with
non-zero intersection number with the ending submanifold $Q$ and hence with
$P$. This contradicts our hypothesis that $P\;$is associated to a canonical
splitting of $G$ over $H$, completing the proof of the lemma.
\end{proof}

\begin{lemma}
\label{canonical Zsplittingimpliescanonicalannulus}An algebraically canonical
splitting of $G=\pi_{1}(M)$ over $H$, which is isomorphic to $\mathbb{Z}$, is
induced by a topologically canonical embedded essential annulus in $M$.
\end{lemma}

\begin{proof}
Lemma \ref{Zendingsubmanifoldforcanonicalsplitting} tells us that if we start
with such a splitting of $G$ over $H$, the associated ending submanifold $P$
of $M_{H}$ can be chosen to have frontier consisting of a single essential
annulus $A$ which carries $H$. We will show that $A$ can be chosen so that its
projection into $M$ is an embedding whose image we again denote by $A$. Thus
$A$ induces the given splitting of $G$ over $H$. As pointed out earlier, it is
trivial that $A$ must be topologically canonical.

In order to prove that we can choose $A$ so as to embed in $M$, it will be
convenient to choose $A$ to be least area in its proper homotopy class. It is
shown in \cite{FHS} that any such least area annulus will also be embedded in
$M_{H}$. (There are two options here. One can choose a Riemannian metric on
$M$ and then choose $A$ to minimise smooth area, see \cite{Meeks-Yau} or
\cite{Nakauchi}, or one can follow the ideas of Jaco and Rubinstein
\cite{JacoRubinstein} and choose a hyperbolic metric on the $2$-skeleton of
some fixed triangulation of $M$ and choose $A$ to minimise $PL$-area. It does
not matter which choice is made.) The fact that $P$ has zero intersection
number with any ending submanifold of any $M_{K}$, with $K\;$isomorphic to
$\mathbb{Z}$, implies, in particular, that the pre-image $X$ of $P$ in
$\widetilde{M}$ does not cross any of its translates. As the frontier $dX$ of
$X$ is an infinite strip covering $A$ and so is connected, it follows that
$dX$ crosses none of its translates. This means that, for all $g\in G$, at
least one of $gdX\cap X$ and $gdX\cap X^{\ast}$ has bounded image in $M_{H}$.
In turn, this implies that the self-intersection number of $A$ is zero so that
$A$ must cover an embedded surface in $M$, by \cite{FHS}. If $A$ properly
covers an embedded annulus, then the stabiliser of $X$ will strictly contain
$H$. But as $P$ is associated to a splitting of $G$ over $H$, we know that the
stabiliser of $X$ is exactly $H$. It follows that either $A$ embeds in $M$ as
required, or that it double covers an embedded $1$-sided Moebius band $B$ in
$M$. In the second case, we can isotop $A$ slightly so that it embeds in $M$
as the boundary of a regular neighbourhood of $B$. Thus in all cases, we can
choose $A$ so that its projection into $M$ is an embedding as required.

This completes the proof that the given algebraically canonical splitting of
$G$ comes from a topologically canonical annulus.
\end{proof}

Next we go through much the same argument with the subgroup $H$ of $G$
replaced by a subgroup $K$ which is isomorphic to $\mathbb{Z}\times\mathbb{Z}%
$. As before, we start with a general discussion of ending submanifolds.
Consider the cover $M_{K}$ of $M$ with $\pi_{1}(M_{K})=K$. The boundary of
$M_{K}$ is incompressible, and so consists of planes, annuli and possibly a
single torus. An embedded incompressible torus $T$ in the interior of $M_{K}$
splits it into two pieces which we will refer to as the left and the right,
and the boundary components of $M_{K}$ are correspondingly split into left and
right boundary components. Let $A_{1}$ and $A_{2}$ be distinct left annulus
components of $\partial M_{K}$. As any loop on a left annulus component of
$\partial M_{K}$ is homotopic into $T$, the Annulus Theorem yields an embedded
incompressible annulus $B_{i}$ joining $A_{i}$ to $T$, $i=1$, $2$. As $A_{1}$
and $A_{2}$ are distinct, it is trivial that the boundary components of
$B_{1}$ and $B_{2}$ which do not lie in $T$ are disjoint. Now it is easy to
show that $B_{1}$ and $B_{2}$ can be isotoped to be disjoint. (Their
intersection must consist of circles and of arcs with both endpoints on $T$.
All nullhomotopic circles can be removed by an isotopy, and then all such arcs
can be removed by an isotopy, starting with innermost arcs. Finally, one can
remove any essential circles by an isotopy.) This implies that all left
annulus components of $\partial M_{K}$ carry the same subgroup of $K$, and the
same comment applies to right annulus components of $\partial M_{K}$. If $P$
is an ending submanifold of $M_{K}$, an annulus component of $dP$ will be
called a left annulus if each of its boundary circles lies in a left annulus
boundary component of $M_{K}$. Right annulus components of $dP$ are defined
similarly. Note that an annulus component of $dP$ may have one boundary circle
in a left annulus component of $\partial M_{K}$ and the other in a right
annulus component, but this can only happen when the left and right annulus
components of $\partial M_{K}$ all carry the same group. We will call such an
annulus \textit{mixed}. It will not separate $M_{K}$, as it can be isotoped to
intersect $T$ in a single essential circle $C$ and then $T$ contains a circle
which meets $C\;$transversely in a single point. But any other properly
embedded compact incompressible annulus $A$ in $M_{K}$ must separate $M_{K}$.
This is because $A$ can be isotoped to be disjoint from $T$, which implies
that one component of $M_{K}-A$ has fundamental group equal to the fundamental
group of $M_{K}$. It can be shown that $M_{K}$ compactifies to $T\times I$,
but we will not need this fact.

\begin{lemma}
\label{ZxZgeneralendingsubmanifold}Any ending submanifold $P$ of $M_{K}$, for
which $K$ is isomorphic to $\mathbb{Z}\times\mathbb{Z}$, is almost equal to
one with the following properties:

\begin{enumerate}
\item The frontier consists only of embedded essential annuli or tori, such
that distinct boundary circles lie in distinct boundary components of $M_{K}$.

\item The frontier has at most one torus component.
\end{enumerate}
\end{lemma}

\begin{proof}
Let $P$ be an ending submanifold of $M_{K}$. As in the proof of Lemma
\ref{Zgeneralendingsubmanifold}, we can arrange that $dP$ consists of
essential annuli and tori, and that distinct boundary circles of $dP$ lie in
distinct boundary components of $M_{K}$. We can also arrange that $dP$ has at
most one torus component, as any two incompressible tori in $M_{K}$ cobound a
product region in $M_{K}$, which can be added to $P$ or subtracted from $P$,
as appropriate. This completes the proof of the lemma.
\end{proof}

\begin{lemma}
\label{ZxZendingsubmanifoldforcanonicalsplitting}Suppose that $P$ is an ending
submanifold of $M_{K}$, for which $K$ is isomorphic to $\mathbb{Z}%
\times\mathbb{Z}$, such that $P$ is associated to a canonical splitting of $G$
over $K$. Then $P$ is almost equal to a connected ending submanifold whose
frontier consists of a single essential torus.
\end{lemma}

\begin{proof}
We can assume that $P$ is chosen as in the immediately preceding Lemma
\ref{ZxZgeneralendingsubmanifold}.

First we show that $dP$ cannot have an annulus component which is mixed. For
suppose that $dP$ has an annulus component joining a left annulus component
$A_{L}$ and a right annulus component $A_{R}$ of $\partial M_{K}$. As
$\partial M_{K}$ has left and right annulus components, it does not have a
torus component so there is an essential torus in $M_{K}$. This torus is the
frontier of an ending submanifold $R$ of $M_{K}$, and we consider the four
intersections $P^{(\ast)}\cap R^{(\ast)}$. None of these can be compact as
each contains a non-compact piece of $A_{L}$ or of $A_{R}$ in its boundary.
This uses the fact that distinct boundary circles of $dP$ lie in distinct
components of $\partial M_{K}$. It follows that $P$ and $R$ have non-zero
intersection number contradicting our hypothesis that $P$ is associated to a
canonical splitting of $G$ over $K$. This contradiction shows that $dP$ does
not have a mixed annulus component. Note that this implies that each component
of $dP$ separates $M_{K}$.

Second we will show that $dP$ cannot have any components which are left or
right annuli. The only remaining possibility will then be that $dP$ consists
of a single torus as claimed. Note that in the following proof, we will only
use the fact that the given splitting is $\mathbb{Z}$-canonical, but not the
fact that it is algebraically canonical. This will allow us to use the same
argument in the last part of the proof of Theorem \ref{specialcanonicaltori}.

Suppose that $dP$ has a left annulus component $A$ which carries the subgroup
$H$ of $K$. Then the pre-image of $A$ in $M_{H}$ consists of infinitely many
copies of $A$. Let $Q$ denote the full pre-image of $P$ in $M_{H}$. Of course,
$Q$ is not an ending submanifold of $M_{H}$ as its frontier $dQ$ is not
compact. The components of $dQ$ may be compact annuli, infinite strips or an
open annulus if $dP$ has a torus component. A compact annulus component of
$dQ$ will be called left or right depending on whether it projects to a left
or right annulus component of $dP$. Recall that distinct boundary circles of
$dP$ lie in distinct boundary components of $M_{K}$. It follows that distinct
boundary circles of $dQ$ lie in distinct boundary components of $M_{H}$. As
$dQ$ contains more than one left annulus component, there is an incompressible
annulus $B$ embedded in $M_{H}$ and joining two such components.

It is possible that $B$ meets other components of $dQ$ in its interior, but we
can isotop $B$ so that it meets each component of $dQ$ in essential circles
only. Hence, if $B$ meets a component $C$ of $dQ$, then $C$ must be an
annulus, compact or open, and must separate $M_{H}$. It follows that if $B$
meets $C$ in more than one essential circle, we can alter $B$ so as to remove
two circles of intersection with $C$, essentially by replacing a sub-annulus
of $B$ by a sub-annulus of $C$. As the boundary circles of $B$ lie in left
annulus components of $dQ$, we can arrange, by repeating this argument, that
$B$ meets only left annulus components of $dQ$. Now $B$ must have a
sub-annulus which joins distinct left annulus components $A_{1}$ and $A_{2}$
of $dQ$ and whose interior does not meet $dQ$. We will replace $B$ by this
sub-annulus and will continue to call this annulus $B$.

Without loss of generality, we can assume that $B$ lies in a component $R$ of
$Q^{\ast}$, and we let $Q_{1}$ and $Q_{2}$ denote the components of $Q$ which
contain $A_{1}$ and $A_{2}$ respectively. This annulus $B$ separates $R$ into
two pieces $U$ and $V$, and the argument of Lemma
\ref{annuluscutsintotwoinfinitepieces} shows that neither $U$ nor $V$ can be
compact. Recall that if $A$ is any left annulus component of $dP$ in $M_{K}$,
then $A$ carries the group $H$ and cuts $M_{K}$ into two pieces one of which
also carries $H$. It follows that any left annulus component of $dQ$ cuts
$M_{H}$ into two pieces one of which projects into $M_{K}$ by a homeomorphism.
As $A_{1}$ and $A_{2}$ are such annuli, it follows that at least one of
$Q_{1}$ and $Q_{2}$ projects into $M_{K}$ by a homeomorphism. In particular,
at least one of $Q_{1}$ and $Q_{2}$ is an ending submanifold of $M_{H}$.
Assume that $Q_{1}$ is an ending submanifold. Cut $M_{H}$ along $A_{1}\cup
B\cup A_{2}$, and let $W$ denote the piece so obtained which contains $X$.
Thus $W$ is also an ending submanifold. If $W^{\ast}$ contains no translate of
$Q_{1}$ other than $Q_{1}$ itself, we remove one such translate from $W$ to
obtain a new ending submanifold $Z$ such that $Z^{\ast}$ contains at least two
translates of $Q_{1}$. Otherwise we let $Z$ equal $W$. Now we let $S=Z\cup
Q_{1}$ and note that $S$ is an ending submanifold of $M_{H}$ and that
$S^{\ast}$ contains a translate of $Q_{1}$. This implies that all four of the
sets $S^{(\ast)}\cap Q^{(\ast)}$ are unbounded. If $X$ and $Y$ denote the
pre-images of $P$ and $S$ respectively in $\widetilde{M}$, it follows that
each of the four sets $X^{(\ast)}\cap Y^{(\ast)}$ has projection into $M_{H}$
which is unbounded, so that $S$ and $P$ have non-zero intersection number.
This contradicts the assumption that $P\;$is associated to a canonical
splitting of $G$ over $K$, so we conclude that $dP$ cannot have a component
which is a left annulus. Similarly, $dP$ cannot have a component which is a
right annulus. This completes the proof of the lemma.
\end{proof}

\begin{lemma}
\label{canonicalZxZsplittingimpliescanonicaltorus}An algebraically canonical
splitting of $G=\pi_{1}(M)$ over $K$, which is isomorphic to $\mathbb{Z}%
\times\mathbb{Z}$, is induced by a topologically canonical embedded essential
torus in $M$.
\end{lemma}

\begin{proof}
Lemma \ref{ZxZendingsubmanifoldforcanonicalsplitting} tells us that if we
start with such a splitting of $G$ over $K$, the associated ending submanifold
$P$ of $M_{K}$ can be chosen to have frontier consisting of a single essential
torus $T$. It is automatic that $T$ carries $K$. We will show that $T$ can be
chosen so that its projection into $M$ is an embedding whose image we again
denote by $T$. Now $T$ induces the given splitting of $G$ over $K$. As pointed
out earlier, it is trivial that $T$ must be topologically canonical.

In order to prove that we can choose $T$ so as to embed in $M$, it will again
be convenient to choose $T$ to be least area. See \cite{Schoen-Yau} or
\cite{JacoRubinstein} for the existence results in the smooth and $PL$ cases
respectively. Again \cite{FHS} shows that a least area torus is embedded. The
fact that $P$ has zero intersection number with any ending submanifold of any
$M_{K}$, with $K\;$isomorphic to $\mathbb{Z}\times\mathbb{Z}$, implies in
particular that the pre-image $X$ of $P$ in $\widetilde{M}$ does not cross any
of its translates. Hence the frontier $dX$ of $X$, which is a plane above $T$,
crosses none of its translates. This means that, for all $g\in G$, at least
one of $gdX\cap X$ and $gdX\cap X^{\ast}$ has bounded image in $M_{K}$. In
turn, this implies that the self-intersection number of $T$ is zero so that
$T$ must cover an embedded surface in $M$, by \cite{FHS}. If $T$ covers an
embedded torus, then the stabiliser of $X$ will contain $K$. But as $P$ is
associated to a splitting of $G$ over $K$, we know that the stabiliser of $X$
is exactly $K$. It follows that either $T$ embeds in $M$ as required, or that
it double covers an embedded $1$-sided Klein bottle in $M$. In this case, we
can isotop $T$ slightly so that it embeds in $M$ as the boundary of a regular
neighbourhood of the Klein bottle. Thus in all cases, we can choose $T$ so
that its projection into $M$ is an embedding as required.

This completes the proof that the given algebraically canonical splitting of
$G$ comes from a topologically canonical torus.
\end{proof}

At this point we have proved part 1 of Theorem
\ref{nonspecialsurfacesandcanonicalsplittings}, by showing that a canonical
splitting of $G$ over $\mathbb{Z}$ or $\mathbb{Z}\times\mathbb{Z}$ comes from
a canonical annulus or torus in $M$. These results were proved in Lemmas
\ref{canonical Zsplittingimpliescanonicalannulus} and
\ref{canonicalZxZsplittingimpliescanonicaltorus}. Now we need to show that a
non-special canonical annulus or torus in $M$ induces a canonical splitting of
$G$.

\begin{lemma}
\label{canonicalannulusimpliescanonicalsplitting}Let $A$ be a topologically
canonical annulus in $M$. Then $A$ determines an algebraically canonical
splitting of $G=\pi_{1}(M)$.
\end{lemma}

\begin{proof}
Let $H$ denote the group carried by $A$, so that we have a splitting of $G$
over $H$. Lift $A$ into $M_{H}$, and let $P$ denote one side of $A$ in $M_{H}%
$. Let $X$ denote the pre-image in $\widetilde{M}$ of $P$. We consider the
annulus boundary components of $M_{H}$. There are two which contain the
boundary circles of $A$. The others lie on the left or right of $A$. As all
these annuli carry exactly the same group, namely $H$, we cannot have such
annuli on both sides of $A$, as this would immediately yield a compact annulus
in $M_{H}$ crossing $A$ in an essential way, contradicting our hypothesis that
$A$ is topologically canonical. Hence we can suppose that $P$ contains no
annulus component of $\partial M_{H}$.

Let $K$ be a subgroup of $G$ isomorphic to $\mathbb{Z}$, and let $Q$ be an
ending submanifold of $M_{K}$ with pre-image $Y$ in $\widetilde{M}$. We will
show that $X$ and $Y$ do not cross. If this is known for every such $K\;$and
$Q$, it will follow that $P$ has intersection number zero with every such $Q$,
as required. Lemma \ref{Zgeneralendingsubmanifold} tells us that we can choose
$Q$ to be a submanifold of $M_{K}$ with frontier consisting only of essential
annuli, and all these annuli carry the same group. Thus $Y$ has frontier
consisting of finitely many infinite strips, which are ``parallel'' in the
sense that each lies within a uniformly bounded distance of each other. We
consider how these strips project into $M_{H}$. If one (and hence every) such
strip has compact image in $M_{H}$, that image will be a (possibly singular)
annulus properly mapped into $M_{H}$. As $P$ contains no annulus component of
$\partial M_{H}$, it follows that each of these annuli can be properly
homotoped to lie in $P^{\ast}$. After this homotopy, one of $Y\cap P$ or
$Y^{\ast}\cap P$ would be empty. It follows that, before the homotopy, one of
these two sets projected to a bounded subset of $M_{H}$, so that $P$ and $Q$
do not cross, as required. Otherwise, the image of each of these strips in
$M_{H}$ is non-compact and the map into $M_{H}$ is proper. Thus each end of
each strip is mapped to the left or to the right side of $A$. If one of these
strips has its two ends on opposite sides of $A$, the corresponding annulus
component of $dQ$ has non-zero intersection number with $A$, contradicting our
hypothesis that $A$ is topologically canonical. It follows that each of these
strips has both of its ends on the same side of $A$, and hence that all the
ends of all the strips are on one side of $A$. Hence they can be properly
homotoped so as to be disjoint from $A$, which again shows that $X$ and $Y$ do
not cross as required.

Now let $K$ be a subgroup of $G$ isomorphic to $\mathbb{Z}\times\mathbb{Z}$,
and let $Q$ be an ending submanifold of $M_{K}$ with pre-image $Y$ in
$\widetilde{M}$. We will show that $X$ and $Y$ do not cross. This will show
that $P\;$has intersection number zero with every such $Q$. Lemma
\ref{ZxZgeneralendingsubmanifold} tells us that we can choose $Q$ so that $dQ$
consists of essential annuli and possibly a torus in $M_{K}$. Thus $Y$ has
frontier $dY$ consisting of infinite strips and possibly a plane. In any case,
if $T$ denotes an incompressible torus in $M_{K}$, it yields a plane $\Pi$ in
$\widetilde{M}$ such that $dY$ lies in a bounded neighbourhood of $\Pi$. Now
we consider how $\Pi$ projects into $M_{H}$. There are two cases. Either $\Pi$
projects into $M_{H}$ to yield a properly immersed plane, or $\Pi$ projects to
a properly immersed annulus. In the first case, as a plane has only one end,
and as a properly immersed plane in $M_{H}$ can only meet $A$ in a compact
set, it follows that the intersection of $\Pi$ with $X$ or $X^{\ast}$ is
compact. As $dY$ lies in a bounded neighbourhood of $\Pi$, it follows that
$dY$ meets $X$ or $X^{\ast}$ in a compact set. This implies that one of the
four intersection sets $X^{(\ast)}\cap Y^{(\ast)}$ is bounded, so that $X$ and
$Y$ do not cross as required. In the second case, there are two subcases
depending on whether the two ends of the annulus image of $\Pi$ are mapped to
the same or opposite sides of $A$. If they are mapped to the same side, then
much as above, it follows that one of the four intersection sets $X^{(\ast
)}\cap Y^{(\ast)}$ has bounded image in $M_{H}$, so that $X$ and $Y$ do not
cross as required. If they are mapped to opposite sides of $A$, it follows
that the torus $T$ in $M_{K}$ has non-zero intersection number with $A$, which
contradicts our hypothesis that $A$ is topologically canonical. This concludes
the proof that a topologically canonical annulus in $M$ determines an
algebraically canonical splitting of $G$.
\end{proof}

Next we prove the corresponding result for a torus. The case of special
canonical tori in $M$ will be discussed at the end of this section.

\begin{lemma}
\label{nonspecialtorusimpliescanonicalsplitting}Let $T$ be a topologically
canonical torus in $M$, which is not special. Then $T$ determines an
algebraically canonical splitting of $G=\pi_{1}(M)$.
\end{lemma}

\begin{proof}
Let $H$ denote the group carried by $T$, so that we have a splitting of $G$
over $H$. Lift $T$ into $M_{H}$, and let $P$ denote one side of $T$ in $M_{H}%
$. Let $X$ denote the pre-image in $\widetilde{M}$ of $P$. As $\partial M_{H}$
cannot have a torus component, it must consist of a (possibly empty)
collection of left annuli and right annuli and planes. If there are any left
annulus components, it follows that the canonical component of $M$ containing
the left side of $T$ is a peripheral Seifert fibre space$.$ Similarly for
right annuli. As we are assuming that $T$ is not special, it follows that
$\partial M_{H}$ cannot have both left annulus and right annulus components.
Without loss of generality, suppose that $\partial M_{H}$ has no left annulus components.

Let $K$ be a subgroup of $G$ isomorphic to $\mathbb{Z}$, and let $Q$ be an
ending submanifold of $M_{K}$ with pre-image $Y$ in $\widetilde{M}$. We will
show that $X$ and $Y$ do not cross. As before, this will show that $P\;$has
intersection number zero with every such $Q$. Lemma
\ref{Zgeneralendingsubmanifold} tells us that we can choose $Q$ to be a
submanifold of $M_{K}$ with frontier consisting only of essential annuli, and
all these annuli carry the same group. Thus the frontier $dY$ of $Y$ consists
of finitely many `parallel' infinite strips. As before there are two cases
depending on whether these strips project properly into $M_{H}$ or project to
compact annuli in $M_{H}$.

Suppose that $dY$ projects to compact annuli in $M_{H}$. The boundary lines of
$dY$ then project to circles in $\partial M_{H}$ which must all lie on the
right side of $T$, as $\partial M_{H}$ has no left annulus components. Thus we
can homotop $dY$ to arrange that $dY$ lies entirely on one side of $T$. It
follows that before the homotopy, one of the four intersection sets
$X^{(\ast)}\cap Y^{(\ast)}$ has bounded image in $M_{H}$, so that $X$ and $Y$
do not cross as required.

Next suppose that $dY$ projects to infinite strips in $M_{H}$. Each boundary
line of each such strip must lie in a boundary component of $M_{H}$. As the
intersection of a strip with $T$ is compact, the two boundary lines of a
single strip must lie on the same side of $T$, so that the strip can be
homotoped to be disjoint from $T$. As the infinite strips forming $dY$ are
parallel, it follows that they all lie on the same side of $T$, and again $X$
and $Y$ do not cross.

Now let $K$ be a subgroup of $G$ isomorphic to $\mathbb{Z}\times\mathbb{Z}$,
and let $Q$ be an ending submanifold of $M_{K}$ with pre-image $Y$ in
$\widetilde{M}$. We will show that $X$ and $Y$ do not cross. As before, this
will show that $P\;$has intersection number zero with every such $Q$. Lemma
\ref{ZxZgeneralendingsubmanifold} tells us that we can assume that $dQ$
consists of essential annuli and possibly one torus. Thus $dY$ consists of
infinite strips and possibly a plane. As in the proof of Lemma
\ref{canonicalannulusimpliescanonicalsplitting}, we consider the projection of
$dY$ into $M_{H}$. An infinite strip component of $dY$ must project to a
properly immersed infinite strip or to an annulus. As above, in either case,
the image of any one strip component of $dY$ can be homotoped to lie on one
side of $T$. If the image of a strip is a singular annulus, its boundary
curves can only lie on the right side of $T$.

Let $S$ denote an incompressible torus in $M_{K}$, and let $\Pi$ denote the
pre-image plane in $\widetilde{M}$. Thus $dY$ lies in a bounded neighbourhood
of $\Pi$. There are three cases for the projection of $\Pi$ into $M_{H}$. We
can obtain a proper map of a plane, torus or annulus.

If $\Pi$ is mapped properly into $M_{H}$, the fact that $\Pi$ has only one end
implies that all except a compact subset of $\Pi$ maps to one side of $T$.
Also the infinite strips in $dY$ must all project to infinite strips each
within a bounded neighbourhood of the image of $\Pi$. It follows that we can
homotop $dY$ to lie on one side of $T$, so that $X$ and $Y$ do not cross.

If the image of $\Pi$ in $M_{H}$ is a singular torus, it follows that every
infinite strip of $dY$ has image a singular annulus and so has boundary lying
on the right side of $T$. Thus $dY\;$can be homotoped into the right side of
$T$, and again $X$ and $Y$ do not cross.

If the image of $\Pi$ in $M_{H}$ is a singular annulus, it will be convenient
to let $\Sigma_{L}$ and $\Sigma_{R}$ denote the canonical pieces of $M$ which
meet the left and right sides of $T$ respectively. Of course, it is possible
that $\Sigma_{L}=\Sigma_{R}$. Now the image of $\Pi$ is a singular annulus
whose two ends must lie on the same side of $T$, as otherwise the tori $S$ and
$T$ would intersect essentially, contradicting our hypothesis that $T$ is a
canonical torus. If an infinite strip in $dY$ has image an infinite strip in
$M_{H}$, that strip must again have its boundary on the same side of $T$ as
the ends of the singular annulus which is the image of $\Pi$. If an infinite
strip in $dY$ has image a singular annulus, we know that its boundary must lie
on the right side of $T$ and so can be homotoped to lie entirely on the right
side of $T$. Thus either all of $dY$ can be homotoped into one side of $T$, so
that $X$ and $Y$ do not cross, or $dY$ has some infinite strips which project
to singular annuli whose boundary curves lie on the right of $T$, but the ends
of the image of $\Pi$ lie on the left of $T$. We will show that this last case
cannot occur. Suppose it does occur. Then $S$ projects to a torus in $M$ which
is homotopic into $\Sigma_{L}$. As the image of $\Pi$ in $M_{H}$ is an
annulus, $S$ and $T$ cannot be homotopic, so it follows that $\Sigma_{L}$ is a
Seifert fibre space. Also there is a singular compact annulus $A$ in $M_{H}$
which has one boundary component on $S$, meets $T$ transversely in a single
circle, and has its second boundary component in an annulus boundary component
of $\partial M_{H}$. The existence of $A$ implies that $\Sigma_{R}$ is a
Seifert fibre space. Further the existence of $A$ shows that $\Sigma_{L}$ and
$\Sigma_{R}$ can be fibred to induce the same fibration on $T$. This
contradicts the fact that $T$ is a canonical torus in $M$, so we deduce that
this last possibility cannot occur.

This completes the proof that a topologically canonical torus in $M$, which is
not special, determines an algebraically canonical splitting of $G$.
\end{proof}

At this point, we have completed the proof of part 2) of Theorem
\ref{nonspecialsurfacesandcanonicalsplittings} by proving Lemmas
\ref{canonicalannulusimpliescanonicalsplitting} and
\ref{nonspecialtorusimpliescanonicalsplitting}. To complete this section, we
need to discuss the algebraic analogue of special canonical tori in $M$. First
we show that these tori never define algebraically canonical splittings.

\begin{lemma}
\label{specialtorusdoesnotimplycanonicalsplitting}If $T$ is a special
canonical torus in $M$, then $\Phi(T)$ is not algebraically canonical.
\end{lemma}

\begin{proof}
Let $T$ be a special canonical torus in $M$, and let $H$ denote the group
carried by $T$, so that we have a splitting of $G$ over $H$. Lift $T$ into
$M_{H}$, and let $P$ denote one side of $T$ in $M_{H}$. Let $X$ denote the
pre-image in $\widetilde{M}$ of $P$. The hypotheses imply that $\partial
M_{H}$ has at least one annulus component on each side of $T$. As before we
will refer to the two sides of $T\;$as the left and right.

Suppose first that $\partial M_{H}$ has at least two annulus components on
each side of $T$. Pick a compact left-annulus $A_{L}$ in $M_{H}$ joining two
distinct left annulus components of $\partial M_{H}$ and let $A_{R}$ be a
similarly defined right annulus. Let $Q_{L}$ denote the component of
$M_{H}-A_{L}$ which does not contain $T$, let $Q_{R}$ denote the component of
$M_{H}-A_{R}$ which does not contain $T$, and let $Q=Q_{L}\cup Q_{R}$. Then
clearly all four of the intersections $P^{(\ast)}\cap Q^{(\ast)}$ are
unbounded, so that $P\;$and $Q$ have non-zero intersection number and so
$\Phi(T)$ is not algebraically canonical in this case.

If $\partial M_{H}$ has only one annulus component on one side of $T$, (or on
both sides of $T)$, we simply make the preceding argument in a double cover of
$M_{H}$ chosen so that each annulus component of $\partial M_{H}$ has two
components above it in the double cover.
\end{proof}

In the above proof, we exhibited an ending submanifold of $M_{H}$ (or of a
double cover of $M_{H}$) which has non-zero intersection number with the
ending submanifold $P$ of $M_{H}$ determined by $T$. This ending submanifold
need not correspond to a splitting of $G$, but the following example shows
that one can often choose it to correspond to a splitting.

\begin{example}
Let $F_{1}$ and $F_{2}$ denote two compact surfaces each with at least two
boundary components. Let $\Sigma_{i}$ denote $F_{i}\times S^{1}$, let $T_{i}$
denote a boundary component of $\Sigma_{i}$, and construct a $3$-manifold $M$
from $\Sigma_{1}$ and $\Sigma_{2}$ by gluing $T_{1}$ to $T_{2}$. Let $T$
denote the torus $\Sigma_{1}\cap\Sigma_{2}$, and let $H$ denote the subgroup
of $G=\pi_{1}(M)$ carried by $T$. Then there is a splitting of $G$ over $H$
which has non-zero intersection number with $\Phi(T)$. If $T_{1}$ is glued to
$T_{2}$ so that the given fibrations by circles do not match, then $T$ is a
special canonical torus in $M$. (If the fibrations are matched, then $M$ is
itself a Seifert fibre space and so has no canonical tori.)

Let $G_{i}$ denote $\pi_{1}(\Sigma_{i})$, and let $C_{i}$ denote the subgroup
of $G_{i}$ carried by $T_{i}$. The starting point of our construction is that
if $F$ is a compact surface with at least two boundary components, and if $S$
denotes a boundary circle of $F$, then $S$ carries an infinite cyclic subgroup
of $\pi_{1}(F)$ which is a free factor of $\pi_{1}(F)$. Now it is easy to give
a splitting of $\pi_{1}(F)$ over $\pi_{1}(S)$, and hence a splitting of
$G_{i}$ over $C_{i}$. If each $\pi_{1}(F_{i})$ is free of rank at least $3$,
then we can write $G_{i}=A_{i}\ast_{C_{i}}B_{i}$. If we let $A$ denote the
subgroup of $G$ generated by $A_{1}$ and $A_{2}$, i.e. $A=A_{1}\ast_{H}A_{2}$,
and define $B$ similarly, then we can express $G$ as $A\ast_{H}B$, and it is
easy to see that this splitting of $G$ has non-zero intersection number with
$\Phi(T)$. If $\pi_{1}(F_{i})$ has rank $2$, then $G_{i}=A_{i}\ast_{C_{i}}$
and a similar construction can be made.
\end{example}

Our next result gives an algebraic characterisation of special canonical tori.

\begin{theorem}
\label{specialcanonicaltori}

\begin{enumerate}
\item If $T$ is a special canonical torus in $M$, carrying the group $H$, then
$\Phi(T)$ is $\mathbb{Z}$-canonical, and $G$ has splittings over two
incommensurable infinite cyclic subgroups of $H$.

\item If $\sigma$ is a $\mathbb{Z}$-canonical splitting of $G$ over a subgroup
$H$ which is isomorphic to $\mathbb{Z}\times\mathbb{Z}$, such that $G$ has
splittings over two incommensurable infinite cyclic subgroups of $H$, then
there is a special canonical torus $T$ in $M$, such that $\Phi(T)=\sigma$.
\end{enumerate}
\end{theorem}

\begin{proof}
1) Lift $T$ into $M_{H}$, and let $P$ denote one side of $T$ in $M_{H}$. Let
$X$ denote the pre-image in $\widetilde{M}$ of $P$. As $T$ is special, the
canonical pieces of $M$ on each side of $T$ are Seifert fibre spaces which
meet $\partial M$ either in tori or in vertical annuli. Hence $\partial M_{H}$
consists of a non-empty collection of left annuli and right annuli and of
planes, with the left annuli and right annuli carrying incommensurable
subgroups of $H$.

Let $K$ be a subgroup of $G$ isomorphic to $\mathbb{Z}$, and let $Q$ be an
ending submanifold of $M_{K}$ with pre-image $Y$ in $\widetilde{M}$. We will
show that $X$ and $Y$ do not cross. As before, this will show that $P\;$has
intersection number zero with every such $Q$. Lemma
\ref{Zgeneralendingsubmanifold} tells us that we can choose $Q$ to be a
submanifold of $M_{K}$ with frontier consisting only of essential annuli, and
all these annuli carry the same group. Thus the frontier $dY$ of $Y$ consists
of finitely many `parallel' infinite strips, i.e. any strip lies in a bounded
neighbourhood of any other strip. As before there are two cases depending on
whether these strips project properly into $M_{H}$ or project to compact
annuli in $M_{H}$.

Suppose that $dY$ projects to compact annuli in $M_{H}$. The boundary lines of
$dY$ then project to circles in $\partial M_{H}$ which must all lie on the
same side of $T$, as the left annuli and right annuli of $\partial M_{H}$
carry incommensurable subgroups of $H$. Thus we can homotop $dY$ to arrange
that $dY$ lies entirely on one side of $T$. It follows that before the
homotopy, one of the four intersection sets $X^{(\ast)}\cap Y^{(\ast)}$ has
bounded image in $M_{H}$, so that $X$ and $Y$ do not cross as required.

Next suppose that $dY$ projects to infinite strips in $M_{H}$. Each boundary
line of each such strip must lie in a boundary component of $M_{H}$. As the
intersection of a strip with $T$ is compact, the two boundary lines of a
single strip must lie on the same side of $T$, so that the strip can be
homotoped to be disjoint from $T$. As the infinite strips forming $dY$ are
parallel, it follows that they all lie on the same side of $T$, and again $X$
and $Y$ do not cross. This completes the proof that $\Phi(T)$ is $\mathbb{Z}$-canonical.

Let $\Sigma_{1}$ and $\Sigma_{2}$ denote the Seifert fibre spaces on each side
of $T$, remembering that it is possible that $\Sigma_{1}$ equals $\Sigma_{2}$.
If $\Sigma_{i}$ has a canonical annulus of $M$ in its frontier, this annulus
will be vertical and so must carry the fibre group of $\Sigma_{i}$. Thus if
each of $\Sigma_{1}$ and $\Sigma_{2}$ are distinct and each has a canonical
annulus of $M$ in its frontier, it is immediate that $G$ has splittings over
two incommensurable infinite cyclic subgroups of $H$. If $\Sigma_{1}$ equals
$\Sigma_{2}$, and has a canonical annulus of $M$ in its frontier, this yields
a splitting of $G$ over an infinite cyclic subgroup of $H$ and a suitable
conjugate of this splitting will be over an incommensurable infinite cyclic
subgroup of $H$. Thus in this case also, $G$ has splittings over two
incommensurable infinite cyclic subgroups of $H$. If $\Sigma_{i}$ has no
canonical annulus in its frontier, it must meet $\partial M$ in a torus
boundary component $T^{\prime}$. In this case, $\Sigma_{i}$ cannot be
homeomorphic to $T\times I$, as this would imply that $T$ was inessential, so
it follows that $\Sigma_{i}$ contains an essential vertical annulus with
boundary in $T^{\prime}$. As before, whether or not $\Sigma_{1}$ and
$\Sigma_{2}$ are distinct, this implies that $G$ has splittings over two
incommensurable infinite cyclic subgroups of $H$, as required. This completes
the proof of part 1) of Theorem \ref{specialcanonicaltori}.

2) Let $\sigma$ be a $\mathbb{Z}$-canonical splitting of $G$ over a subgroup
$H$ which is isomorphic to $\mathbb{Z}\times\mathbb{Z}$, such that $G$ has
splittings over two incommensurable infinite cyclic subgroups of $H$. As
usual, we consider the cover $M_{H}$ of $M$, whose boundary must consist of at
most one torus, some planes and some left annuli and right annuli. If $G$
admits a splitting over an infinite cyclic subgroup $L$, it follows that there
is an embedded essential annulus $A$ in $M$ such that $A$ carries a subgroup
of $L$. Thus the fact that $G$ admits splittings over two incommensurable
infinite cyclic subgroups of $H$, implies that $M$ contains essential annuli
$A$ and $B$ which also carry incommensurable subgroups of $H$. In particular,
they must lift to annuli in $M_{H}$. If $M_{H}$ has a torus boundary component
$T$, then any loop in $\partial M_{H}$ is homotopic into $T$. It follows that
$M_{H}$ contains essential annuli $A^{\prime}$ and $B^{\prime}$ (which may not
be lifts of $A$ and $B)$ which carry incommensurable subgroups of $H$ and each
have one boundary component in $T$. But this implies that $M_{H}$ is
homeomorphic to $T\times I$, and hence that $H\;$has finite index in $G$,
which contradicts the hypothesis that $G$ splits over $H$. This contradiction
shows that $M_{H}$ does not have a torus boundary component. It follows that
$\partial M_{H}$ has both left annuli and right annuli, and the left and right
annuli of $\partial M_{H}$ carry incommensurable subgroups of $H$.

Let $P$ denote an ending submanifold for the given $\mathbb{Z}$-canonical
splitting $\sigma$ of $G$, where $P$ is chosen as in Lemma
\ref{ZxZgeneralendingsubmanifold}. As the left and right annulus components of
$\partial M_{H}$ carry incommensurable subgroups of $H$, it follows that no
component of $dP$ can be a mixed annulus. Now we use the same argument as in
the second part of the proof of Lemma
\ref{ZxZendingsubmanifoldforcanonicalsplitting}\ to show that as $\sigma$ is
$\mathbb{Z}$-canonical, $dP$ cannot have a component which is a left annulus
nor a right annulus. It follows that $dP$ consists of a single essential torus
$T$ in $M_{H}$, and the argument of Lemma
\ref{canonicalZxZsplittingimpliescanonicaltorus} shows that $T$ can be chosen
to project to an embedding in $M$. Finally, the fact that $\partial M_{H}$ has
both left annuli and right annuli implies that $T$ is a special canonical
torus as required. This completes the proof of part 2) of Theorem
\ref{specialcanonicaltori}.
\end{proof}

We can summarise the results of this section as follows.

\begin{theorem}
\label{bijectionsummary}Let $M$ be an orientable Haken $3$-manifold with
incompressible boundary, and let $\Phi$ denote the natural map from the set of
isotopy classes of embedded essential annuli and tori in $M$ to the set of
splittings of $G=\pi_{1}(M)$ over a subgroup isomorphic to $\mathbb{Z}$ or to
$\mathbb{Z}\times\mathbb{Z}$. Then

\begin{enumerate}
\item $\Phi$ induces a bijection between the non-special canonical annuli and
tori in $M$ and the canonical splittings of $G$, and

\item $\Phi$ induces a bijection between the special canonical tori in $M$ and
those $\mathbb{Z}$-canonical splittings of $G$ over a subgroup $H$ which is
isomorphic to $\mathbb{Z}\times\mathbb{Z}$, such that $G$ has splittings over
two incommensurable infinite cyclic subgroups of $H$.
\end{enumerate}
\end{theorem}

\section{The Deformation Theorem}

\label{deformationtheoremchapter}

In this section, we will give our proof of Johannson's Deformation Theorem
based on the preceding work in this paper. We will give a brief comparison of
our arguments with the previous proofs after our proof of Theorem
\ref{deformationtheorem}.

The results so far imply the following :

\begin{theorem}
\label{samegraphofgroups}Let $M$ and $N$ be orientable Haken $3$-manifolds
with incompressible boundary, and let $\Gamma_{M}$ and $\Gamma_{N}$ denote the
graphs of groups structures for $\pi_{1}(M)$ and $\pi_{1}(N)$ determined by
the canonical annuli and tori in $M$ and $N$ respectively. If $f:M\rightarrow
N$ is a homotopy equivalence, the induced isomorphism $f_{\ast}:\pi
_{1}(M)\rightarrow\pi_{1}(N)$ yields a graph of groups isomorphism $\Gamma
_{M}\rightarrow\Gamma_{N}$.
\end{theorem}

\begin{proof}
Let $G$ denote $\pi_{1}(M)$. Theorem \ref{bijectionsummary} tells us that
there is a natural bijection between the collection of isotopy classes of all
canonical annuli and tori in $M$ and a certain collection of conjugacy classes
of splittings of $G$, some canonical and some $\mathbb{Z}$-canonical. As this
collection of splittings of $G$ is defined purely algebraically, it follows
that the conjugacy classes of splittings of $G$ obtained from $M$ and from $N$
are the same. Now Theorem
\ref{splittingsarecompatibleiffzerointersectionnumber} implies that they yield
isomorphic graphs of groups structures for $G$, as required.
\end{proof}

Before we start discussing our proof of the Deformation Theorem, we need to
consider more carefully the pieces obtained by splitting $M$ along its
JSJ-system $\mathcal{F}(M)$. These pieces $M_{j}$ may have compressible
boundary. However, if we denote $M_{j}\cap\partial M$ by $\partial_{0}M_{j}$
and denote the frontier $frM_{j}$ by $\partial_{1}M_{j}$, then each of
$\partial_{0}M_{j}$ and $\partial_{1}M_{j}$ is incompressible in $M$ and
$\partial M_{j}=\partial_{0}M_{j}\cup\partial_{1}M_{j}$. In what follows, we
will consider a triple $(L,\partial_{0}L,\partial_{1}L)$, where $L$ is a
compact orientable irreducible $3$-manifold whose boundary is divided into
subsurfaces $\partial_{0}L$ and $\partial_{1}L$ (which need not be connected)
such that $\partial_{0}L$ is incompressible in $L$. Often $\partial_{1}L$ will
also be incompressible. We will say that the pair $(L,\partial_{0}L)$ is
\textit{simple} if every embedded incompressible annulus or torus in
$(L,\partial_{0}L)$ is parallel into $\partial_{0}L$ or into $\partial_{1}L$.
This is the commonly used definition but it is different from the one used by
Johannson in \cite{J}. (His definition is in terms of characteristic
submanifolds of Haken manifolds with boundary patterns.) We will say that
$(L,\partial_{0}L)$ is \textit{fibred} if either $L$ is an $I$-bundle over a
surface $F$ such that $\partial_{1}L$ is the restriction of this $I$-bundle to
$\partial F$, or if $L$ is Seifert fibred and $\partial_{0}L$ consists of tori
and vertical annuli in $\partial L$.

We will want to use the JSJ-system $\mathcal{F}(M)$ of $M$ rather than the
characteristic submanifold $V(M)$. Thus the following characterisation of
$\mathcal{F}(M)$ will be very useful.

\begin{theorem}
Let $M$ denote an orientable Haken $3$-manifold with incompressible boundary,
and let ${S_{1},...,S_{m}}$ be a family of disjoint essential annuli and tori
embedded in $M$. Let $(M_{j},\partial_{0}M_{j},\partial_{1}M_{j})$ denote the
manifolds obtained by splitting $M$ along the $S_{i}$. This family is the
JSJ-system $\mathcal{F}(M)$ of $M$ if and only if

\begin{enumerate}
\item each $(M_{j},\partial_{0}M_{j})$ is either simple or fibred, and

\item the system ${S_{1},...,S_{m}}$ is minimal with respect to the property 1).
\end{enumerate}
\end{theorem}

This theorem is proved towards the end of \cite{JS}. Their development uses
hierarchies, and a definition of the characteristic submanifold which is
different from the above characterisation, and they prove the Annulus, Torus
and Enclosing Theorems in the process. In \cite{Scott:Strong}, Scott takes the
above properties to define the JSJ-system, proves strong versions of the
Annulus and Torus Theorems and then deduces the Enclosing Theorem.

Recall that in Proposition \ref{JSJ-systemequalsF}, we described how to obtain
$\mathcal{F}(M)$ from $V(M)$. Now we will describe how to obtain $V(M)$ from
the JSJ-system $\mathcal{F}(M)$. First we will construct a new family
$\mathcal{F}^{\prime}$ of essential annuli and tori in $M$ which will be the
frontier of the characteristic submanifold $V(M)$. In order to see the need
for this, suppose that $M$ consists of two Seifert fibre spaces glued together
along a boundary torus $T$ so that their fibrations do not match. Then
$\mathcal{F}(M)$ consists of $T$, but $V(M)$ has two components homeomorphic
to the two constituent Seifert fibre spaces of $M$. Thus the frontier of
$V(M)$ consists of two copies of $T$, and the two components of $V(M)\;$are
separated by the product region between these two tori. Also if $M$ consists
of two annulus free hyperbolic manifolds glued together along a boundary torus
$T$, then $T$ is the only essential torus in $M$ and $M$ has no essential
annuli. Thus $\mathcal{F}(M)$ consists of $T$, and $V(M)$ consists of a
product neighbourhood of $T$. Again the frontier of $V(M)$ consists of two
copies of $T$. Similar examples occur when one glues manifolds along an annulus.

We construct $\mathcal{F}^{\prime}$ as follows. If both sides of $S_{k}$ are
fibred, we add a second parallel copy of $S_{k}$, thus creating an extra
complementary component $X_{k}$, which is homeomorphic to $S_{k}\times I$ and
$\partial_{1}X_{k}$ consists of these two copies of $S_{k}$. If both sides of
$S_{k}$ are simple and not fibred, we also add a second parallel copy of
$S_{k}$. Splitting $M$ along $\mathcal{F}^{\prime}$ yields a manifold which is
homeomorphic to the result of splitting $M$ along $\mathcal{F}$ apart from
some extra components $X_{k}$. It remains to assign the components of $M$
split along $\mathcal{F}^{\prime}$ to be part of $V(M)$ or part of
$\overline{M-V(M)}$. If both sides of $S_{k}$ were fibred, we assign $X_{k}$
to $\overline{M-V(M)}$. If both sides of $S_{k}$ were simple and not fibred,
we assign $X_{k}$ to $V(M)$. The other components are assigned to $V(M)$ if
they are fibred and to $\overline{M-V(M)}$ otherwise. It is clear that no two
adjacent pieces are assigned to $V(M)$ or to $\overline{M-V(M)}$. It follows
that if a piece $(M_{j},\partial_{0}M_{j})$ of $\overline{M-V(M)}$ is fibred,
then it must be one of the extra pieces $X_{k}$ described above. Thus
$(M_{j},\partial_{0}M_{j})$ is homeomorphic to $(S^{1}\times S^{1}\times
I,\emptyset)$ or to $(S^{1}\times I\times I,S^{1}\times\partial I\times I)$.

We will need to use the following important result of Waldhausen.

\begin{theorem}
(Waldhausen's Homeomorphism Theorem) If $f:(M,\partial M)\rightarrow
(N,\partial N)$ is a homotopy equivalence of Haken manifolds which induces a
homeomorphism on $\partial M$, then $f$ can be deformed relative to $\partial
M$ to a homeomorphism.
\end{theorem}

We will note the following extension due independently to Evans \cite{Evans},
Swarup \cite{Sw} and Tucker \cite{Tucker}. This will only be used to clear up
some loose ends.

\begin{theorem}
\label{EvansSwarupTucker}If $f:M\rightarrow N$ is a proper map of Haken
manifolds inducing an injection of fundamental groups, then there is a proper
homotopy $f_{t}:M\rightarrow N$, with $f_{0}=f$, such that one of the
following holds:

\begin{enumerate}
\item $f_{1}$ is a covering map,

\item $M$ is a $I$-bundle over a closed surface and $f_{1}(M)\subset\partial
N$,

\item $N$ (and hence also $M$) is a solid torus or a solid Klein bottle and
$f_{1}$ is a branched covering with branch set a circle,

\item $M$ is a (possibly non-orientable) handlebody and $f_{1}(M)\subset
\partial N$.
\end{enumerate}

If $B$ and $C$ are components of $\partial M$ and $\partial N$ respectively
such that $f\mid B:B\rightarrow C$ is a covering map, then we can choose
$f_{t}\mid B=f\mid B$, for all $t$.
\end{theorem}

Now we can state and prove Johannson's Deformation Theorem.

\begin{theorem}
\label{deformationtheorem}(Johannson's Deformation Theorem) Let $M$ and $N$ be
orientable Haken $3$-manifolds with incompressible boundary, and let $V(M)$
denote the characteristic submanifold of $M$. If $f:M\rightarrow N$ is a
homotopy equivalence, then $f$ is homotopic to a map $g$ such that
$g(V(M))\subset V(N)$, the restriction of $g$ to $V(M)$ is a homotopy
equivalence from $V(M)$ to $V(N)$, and the restriction of $g$ to
$\overline{M-V(M)}$ is a homeomorphism onto $\overline{N-V(N)}$.
\end{theorem}

\begin{proof}
We start by applying Theorem \ref{samegraphofgroups}. Let $\mathcal{F}%
(M)\;$denote the JSJ-system of $M$. It follows immediately from Theorem
\ref{samegraphofgroups} that we can homotop $f$ to a map $g$ such that $g$
maps $\mathcal{F}(M)$ to $\mathcal{F}(N)$ by a homeomorphism and
$g^{-1}\mathcal{F}(N)=\mathcal{F}(M)$. Further, it is now automatic that $g$
induces a bijection between the components of $M-\mathcal{F}(M)$ and of
$N-\mathcal{F}(N)$, and that if $X$ is the closure of a component of
$M-\mathcal{F}(M)$, and if $Y$ denotes the closure of the component of
$N-\mathcal{F}(N)$ which contains $g(X)$, then $g\mid X:X\rightarrow Y$ is a
homotopy equivalence. If $X$ does not meet $\partial M$, then $g\mid
X:X\rightarrow Y$ is proper, i.e. $g(\partial X)\subset\partial Y$, and so
Waldhausen's Homeomorphism Theorem implies that it is homotopic to a
homeomorphism by a homotopy fixed on $\partial X$. Note that this is true
whether or not $X$ is fibred, so that the restriction of $g$ to interior
components of $V(M)$ must also be a homeomorphism.

In order to complete the proof of the Deformation Theorem, it remains to
handle those components of $M-\mathcal{F}(M)$ which are not fibred and meet
$\partial M$. Let $M_{j}$ be the closure of such a component, and let $N_{j}$
denote the closure of the component of $N-\mathcal{F}(N)$ which contains
$g(M_{j})$. Let $g_{j}$ denote the induced map $M_{j}\rightarrow N_{j}$ and
recall that $g_{j}$ restricts to a homeomorphism of $\partial_{1}M_{j}$ with
$\partial_{1}N_{j}$. It was shown in \cite{NS} that any simple pair
$(M_{j},\partial_{0}M_{j})$ such that $(M_{j},\partial_{1}M_{j})$ contains an
embedded $\pi_{1}$-injective annulus not parallel into $\partial_{1}M_{j}$ is
fibred. In particular, it follows that no component of $\partial_{0}M_{j}$ is
an annulus. Now Theorem \ref{Jacoresult}, which we prove below, tells us that
$g_{j}$ can be homotoped to a homeomorphism while keeping $\partial_{1}M_{j}$
mapped into $\partial_{1}N_{j}$ during the homotopy. Note that the final
homeomorphism may flip the $I$-factor of some of the annuli in $\partial
_{1}M_{j}$.
\end{proof}

At this point, we can give the promised comparison of the proofs of the
Deformation Theorem. Johannson proved his result by directly considering the
given homotopy equivalence $f:M\rightarrow N$ and describing a sequence of
homotopies whose end result is $g$. In \cite{Swarup:Deformation}, Swarup gave
a very simple proof of the special case of the Deformation Theorem when $M$ is
annulus free. (See also Remark \ref{Swarupresult} below.) This was an
interesting application of the theory of ends. In \cite{JacoCBMSNotes}, Jaco
gave an alternative proof of the Deformation Theorem which was stimulated by
Swarup's ideas. But he was unable to find a direct generalisation of Swarup's
result on ends. As in the argument which we outlined above, but using
completely different methods, Jaco's first step was to homotop $f$ to a map
$g$ such that $g$ maps $\mathcal{F}(M)$ to $\mathcal{F}(N)$ by a homeomorphism
and $g^{-1}\mathcal{F}(N)=\mathcal{F}(M)$, and his second step was to deal
with the non-fibred pieces of $M$ which meet $\partial M$. A special case of
Theorem \ref{Jacoresult}, sufficient for this second step, was proved by Jaco
in \cite{JacoCBMSNotes}, but our proof is rather different, and we correct
some omissions in his argument.

Now we return to the last part of the proof of Theorem
\ref{deformationtheorem} which handles those components of $M-\mathcal{F}(M)$
which are not fibred and meet $\partial M$. The first step in this argument is
to deform $\partial_{0}M_{j}$ into $\partial N_{j}$. We prove the following
which is far more general than is required for this particular problem, but
seems to be of independent interest.

\begin{lemma}
\label{generalisedJacolemma}Let $(M,\partial_{0}M,\partial_{1}M)$,
$(N,\partial_{0}N,\partial_{1}N)$ be orientable $3$-manifolds with $M$ and $N$
Haken and $\partial_{0}M$ incompressible. Let $f:M\rightarrow N$ be any
homotopy equivalence such that the restriction $f\mid\partial_{1}M$ is a
homeomorphism onto $\partial_{1}N$. Let $t$ be a component of $\partial_{0}M$
and suppose that no $\pi_{1}$-injective map of an annulus into $(M,\partial
_{0}M)$ is essential in $(M,\partial M)$. Then $f\mid t$ can be deformed
relative to $\partial t$ into $\partial N$.
\end{lemma}

\begin{remark}
\label{Swarupresult}In \cite{Swarup:Deformation}, Swarup proved a result which
is the special case of this lemma when $\partial_{1}M$ is empty and so $t$ is
closed. He used it to give a simple proof of Johannson's Deformation Theorem
in the special case when $M$ does not admit any essential annulus. It has been
a problem for many years to find the `correct' generalisation of Swarup's
result to the case of surfaces with boundary. This lemma seems to be that generalisation.
\end{remark}

\begin{proof}
Let $p:M_{t}\rightarrow M$ be the cover of $M$ corresponding to the image of
$\pi_{1}(t)$ in $\pi_{1}(M)$ and $q:N_{t}\rightarrow N$ the cover of $N$
corresponding to the image of $f_{\ast}\pi_{1}(t)$. We have $f_{t}%
:M_{t}\rightarrow N_{t}$, a lift of $f$. We denote by $\widetilde{t}$, a lift
of $t$ at the base point of $M_{t}$, so that $\widetilde{t}$ is homeomorphic
to $t$, and observe that it is enough to deform $f_{t}\mid\widetilde{t}$ into
$\partial N_{t}$. We denote by $\partial_{0}M_{t}$ the inverse image of
$\partial_{0}M$ under $p$. We next consider the map induced by inclusion on
homology groups with integer coefficients:
\[
i_{\ast}:H_{2}(M_{t},\partial_{0}M_{t})\rightarrow H_{2}(M_{t},\partial
M_{t}).
\]
We claim that $i_{\ast}$ is the zero map. If $\partial_{1}M_{t}$ is empty, the
above assertion is equivalent to showing that $H_{2}(M_{t},\partial M_{t})=0$,
the case that was considered in \cite{Swarup:Deformation}. In the general
case, the proof proceeds in a similar fashion. Consider a component $s$ of
$\partial_{0}M_{t}$ other than $\widetilde{t}$. If $\pi_{1}(s)$ is
non-trivial, we have a homotopy annulus joining any loop in $s$ to a loop in
$\widetilde{t}$. By our assumption any such annulus can be homotoped relative
to its boundary into $\partial M_{t}$. It follows that any loop in $s$ is
peripheral in $s$ and thus $s$ is either $S^{1}\times I$ or $S^{1}\times I$
with a closed subset of $S^{1}\times{1}$ removed. Moreover, there is an
annulus in $\partial M_{t}$ joining any non-trivial loop in $s$ to a boundary
component of $\widetilde{t}$.

We next observe that $H_{2}(M_{t},\partial_{0}M_{t})$ is generated by properly
embedded, two sided surfaces in $(M_{t},\partial_{0}M_{t})$. If $k$ is such a
surface, consider a component $l$ of $\partial k$ which is not in
$\widetilde{t}$. Let $s$ denote the component of $\partial_{0}M_{t}$ which
contains $l$. If $l$ is contractible in $s$, it must bound a $2$-disc in $s$
which we add to $k$ and then push off $s$. If $l$ is essential in $s$, we add
to $k$ an annulus $a$ in $\partial M$ which joins $l$ to a peripheral loop
$l^{\prime}$ in $\widetilde{t}$ and push $a-l^{\prime}$ off $\partial M_{t}$.
Thus, in either case, we can modify $k$ without changing its image in
$H_{2}(M_{t},\partial M_{t})$ so as to remove the component $l$. By repeating
this process, we see that the image of $H_{2}(M_{t},\partial_{0}M_{t})$ in
$H_{2}(M_{t},\partial M_{t})$ is generated by surfaces $k$ whose boundary is
in $\widetilde{t}$. Now any such surface is homologous to an incompressible
surface, and we observe that any incompressible surface in $M_{t}$ with
boundary in $\widetilde{t}$ is parallel into $\widetilde{t}$. Thus the image
of $H_{2}(M_{t},\partial_{0}M_{t})$ in $H_{2}(M_{t},\partial M_{t})$ is zero.

Using Poincar\'{e} duality for cohomology with compact supports, it follows
that the induced map in cohomology from $H_{c}^{1}(M_{t},\partial_{1}M_{t})$
to $H_{c}^{1}(M_{t})$ is zero. As $f$ is a homotopy equivalence, the induced
map $f_{t}:M_{t}\rightarrow N_{t}$ is a proper homotopy equivalence, so that
it induces an isomorphism $H_{c}^{1}(N_{t})\rightarrow H_{c}^{1}(M_{t})$.
Further the hypotheses imply that $f_{t}$ induces a homeomorphism on
$\partial_{1}M_{t}$, so that it also induces an isomorphism $H_{c}^{1}%
(N_{t},\partial_{1}N_{t})\rightarrow H_{c}^{1}(M_{t},\partial_{1}M_{t})$.
Applying Poincar\'{e} duality in $N_{t}$, we see that the induced map in
homology (with integer coefficients):
\[
j_{\ast}:H_{2}(N_{t},\partial_{0}N_{t})\rightarrow H_{2}(N_{t},\partial
N_{t})
\]
is also the zero map. Next consider $f_{t}(\widetilde{t})$. Since $f_{t}$
induces an isomorphism $\pi_{1}(\widetilde{t})\rightarrow\pi_{1}(N_{t})$ and
is already an embedding on $\partial\widetilde{t}$, we can homotop $f_{t}%
\mid\widetilde{t}$ relative to $\partial\widetilde{t}$ to become an embedding.
One method to prove this is to choose a metric on $N_{t}$ which blows up away
from $f_{t}(\widetilde{t})$, and homotop $f_{t}\mid\widetilde{t}$ to be of
least area. A more classical argument involves taking a regular neighborhood
$W$ of $f_{t}(\widetilde{t})$, and compressing its frontier $frW$ by adding
$2$-handles or by cutting $W$ along a compressing disc for $frW$. One can
homotop $f_{t}\mid\widetilde{t}$ so as to avoid these compressing discs, so
that one obtains a submanifold $W^{\prime}$ of $N_{t}$ which contains
$f_{t}(\widetilde{t})$ and has incompressible frontier. We let $t^{\prime}$
denote a component of $frW^{\prime}$. As $f_{t}$ induces an isomorphism
$\pi_{1}(\widetilde{t})\rightarrow\pi_{1}(N_{t})$, it follows that there is a
homotopy inverse map $N_{t}\rightarrow\widetilde{t}$. This map will induce a
proper map $t^{\prime}\rightarrow\widetilde{t}$ which induces an injection of
fundamental groups and an injection of boundaries. It follows that this map is
homotopic to a homeomorphism rel $\partial t^{\prime}$. Hence $f_{t}%
\mid\widetilde{t}$ can be homotoped to an embedding with image $t^{\prime}$ as
required. Since $f_{t}(\widetilde{t})$, which is now embedded, represents the
zero element in $H_{2}(N_{t},\partial N_{t})$, we see that $f_{t}%
\mid\widetilde{t}$ is parallel into $\partial N_{t}$. Thus $f_{t}%
\mid\widetilde{t}$ can be deformed into $\partial N_{t}$ relative to
$\partial\widetilde{t}$ and hence $f\mid t$ can be deformed into $\partial N$
relative to $\partial t$. This completes the proof of Lemma
\ref{generalisedJacolemma}.
\end{proof}

\begin{remark}
The condition that $i^{\ast}:H_{c}^{1}(M_{t},\partial_{1}M_{t})\rightarrow
H_{c}^{1}(M_{t})$ is zero can be formulated in algebraic terms. In the case
when $\partial_{1}M_{t}=\emptyset$, the condition becomes that $H_{c}%
^{1}(M_{t})=0$ which is equivalent to $H^{1}(G,\mathbb{Z}[G/H])=0$, where
$G=\pi_{1}(M)$ and $H$ is the image of $\pi_{1}(t)$ in $G=\pi_{1}(M)$. In the
general case, the condition is that the induced map:
\[
H^{1}(G,H_{1},...,H_{m};\mathbb{Z}[G/H])\rightarrow H^{1}(G,\mathbb{Z}[G/H])
\]
is zero,where $H_{i}$ represent the subgroups of $G$ corresponding to the
components of $\partial_{1}M$. The condition $H^{1}(G,\mathbb{Z}[G/H])=0$
implies that the number of ends of the pair $(G,H)$ is one. This last
implication was observed by Jaco in \cite{JacoCBMSNotes}.
\end{remark}

Before continuing we need to know more about when the existence of a singular
essential annulus in a $3$-manifold implies the existence of an embedded one.
The following, which is the Relative Annulus Theorem in \cite{Scott:Strong},
provides the information we need.

\begin{theorem}
\label{relativeannulustheorem}(Relative Annulus Theorem) Let $(M,\partial
_{0}M)$ be an orientable Haken manifold with each component of $\partial_{0}M$
incompressible. Let $\alpha:(A,\partial A)\rightarrow(M,\partial_{0}M)$ be a
$\pi_{1}$-injective map of an annulus which cannot be homotoped relative to
$\partial A$ into $\partial M$. Then, either

\begin{enumerate}
\item there is an embedded annulus in $(M,\partial_{0}M)$ which is not
parallel into $\partial M$, or

\item $\partial_{0}M$ is the disjoint union of some annuli and a surface $T$
such that $(M,T)$ is an $I$-bundle. Moreover, the only $(M,T)$ in this case
which do not satisfy 1) are the product bundle over the twice punctured disc
and the twisted $I$-bundle over the once punctured Mobius band.
\end{enumerate}
\end{theorem}

Now we are ready to prove Theorem \ref{Jacoresult} which we used at the end of
our proof of Johannson's Deformation Theorem. Jaco stated a special case of
this result as Lemma X.23 of \cite{JacoCBMSNotes}. The special case was that
he assumed that $\partial_{1}M$ consists of annuli and tori. Of course, this
is sufficient for the application to prove the Deformation Theorem. However,
there are some omissions in his proof. The main problem is that the statement
of the theorem has the assumption that any embedded incompressible annulus in
$(M,\partial_{0}M)$ is parallel into $\partial_{0}M$ or $\partial_{1}M$,
whereas in his proof he assumes that any $\pi_{1}$-injective map of the
annulus into $(M,\partial_{0}M)$ is homotopic into $\partial_{0}M$ or
$\partial_{1}M$. These two conditions are not equivalent, but Theorem
\ref{relativeannulustheorem} above tells us that there are a very small number
of exceptional cases. We give our own proof of Jaco's result below. In the
non-exceptional cases of Lemma \ref{relativeannulustheorem}, we replace the
first part of Jaco's argument by the result of Lemma
\ref{generalisedJacolemma}. In the exceptional cases, neither his argument nor
that result apply. We give a special argument for these cases.

\begin{theorem}
\label{Jacoresult}Let $(M,\partial_{0}M,\partial_{1}M)$, $(N,\partial
_{0}N,\partial_{1}N)$ be orientable $3$-manifolds with $M$ and $N$ Haken, and
$\partial_{0}M$ incompressible. Suppose that no component of $\partial_{0}M$
or $\partial_{0}N$ is an annulus and that any embedded incompressible annulus
in $(M,\partial_{0}M)$ is parallel into $\partial_{0}M$ or $\partial_{1}M$. If
$f:M\rightarrow N$ is a homotopy equivalence such that the restriction
$f\mid\partial_{1}M$ is a homeomorphism onto $\partial_{1}N$, then $f$ can be
deformed to a homeomorphism while keeping $\partial_{1}M$ mapped into
$\partial_{1}N$ during the homotopy.
\end{theorem}

\begin{proof}
We start by considering a component $t$ of $\partial_{0}M$, and observing that
either Lemma \ref{generalisedJacolemma} applies, or we are in the special case
2) in Theorem \ref{relativeannulustheorem}.

First suppose that the special case of Theorem \ref{relativeannulustheorem}
does not arise, so that any $\pi_{1}$-injective annulus in $(M,\partial_{0}M)$
is properly homotopic into $\partial M$. Then Lemma \ref{generalisedJacolemma}
applies to every component of $\partial_{0}M$, so we can homotop $f$ rel
$\partial_{1}M$ to a proper map. We will continue to call this map $f$.

Next we claim that, for each component $t$ of $\partial_{0}M$, we can deform
$f\mid t$ to an embedding so that during the homotopy each component of
$\partial t$ either stays fixed or is moved across an annulus component of
$\partial_{1}N$. The argument here uses the fact that $f\mid t$ is $\pi_{1}%
$-injective, but this is not enough. It is also necessary to use the fact that
any $\pi_{1}$-injective annulus in $(M,\partial_{0}M)$ is properly homotopic
into $\partial M$. This is another omission in Jaco's proof. We let $\Sigma$
denote the component of $\partial N$ which contains $f(t)$.

If $t$ is closed, then $f\mid t$ is homotopic to a covering map of $\Sigma$ of
some degree $d$. Theorem 1.3 of \cite{Scott:HakenM3} tells us that, as
$M\;$and $N$ are orientable, $d$ must equal $1$, so that we have deformed
$f\mid t$ to a homeomorphism onto its image, as required.

Next we consider the case when $t$ has non-empty boundary. We already know
that $f\mid t$ embeds $\partial t$ in $\Sigma$. Now we consider the lift
$f_{t}$ of $f\mid t$ into the cover $\Sigma_{t}$ of $\Sigma$ whose fundamental
group equals $f_{\ast}\pi_{1}(t)$. We can homotop $f_{t}$ rel $\partial t$ to
an embedding. This homotopy induces a homotopy of $f\mid t$. We will now
assume that this homotopy has been done and use the same notation for the new
maps. Thus $f_{t}$ is an embedding in $\Sigma_{t}$, and we denote its image by
$X$. Let $\widetilde{\Sigma}$ denote the universal cover of $\Sigma$, and let
$\widetilde{X}$ denote the pre-image in $\widetilde{\Sigma}$ of $X$. The full
pre-image in $\widetilde{\Sigma}$ of $t$ consists of $\widetilde{X}$ and all
its translates by $\pi_{1}(\Sigma)$. Suppose that $f\mid t$ is not an
embedding. Then some translate $g\widetilde{X}$ of $\widetilde{X}$ must meet
$\widetilde{X}$ but not equal $\widetilde{X}$. Let $\widetilde{Y}$ denote a
component of this intersection, and let $Y$ denote the image of $\widetilde
{Y}$ in $\Sigma_{t}$. This will be a compact subsurface of $X$. Note that the
boundary of $Y$ consists of circles which project to $\partial t$. Consider an
essential (possibly singular) loop $C$ in $Y$ and let $l$ denote a line in
$\widetilde{Y}$ above $C$. Let $\alpha$ denote a generator of the stabiliser
of $l$. Then $\alpha$ lies in the stabiliser $H$ of $\widetilde{X}$ and in the
stabiliser $gHg^{-1}$ of $g\widetilde{X}$. Thus there is $\beta$ in $H$ such
that $\alpha=g\beta g^{-1}$. As $g\widetilde{X}$ is not equal to
$\widetilde{X}$, we know that $g$ does not lie in $H$. Thus we obtain a
$\pi_{1}$-injective (possibly singular) annulus $A$ in $M$ with both ends on
$t$ and one end at $C$, such that $A$ is not properly homotopic into $t$. But
any $\pi_{1}$-injective annulus in $(M,\partial_{0}M)$ is properly homotopic
into $\partial M$. It follows that some component $D$ of the closure of
$\partial M-t$ is an annulus and that $A$ can be homotoped to cover $D$,
keeping $\partial A$ in $t$ during the homotopy. As $D$ is a union of
components of $\partial_{0}M$ and $\partial_{1}M$, and as no component of
$\partial_{0}M$ is an annulus, it follows that $D$ is a component of
$\partial_{1}M$. We let $E$ denote the component $f(D)$ of $\partial_{1}N$. A
particular consequence of the preceding argument is that $C$ must be homotopic
in $t$ into a component of $\partial t$. Hence every essential loop on $Y$ is
homotopic into a component of $\partial X$, so that $Y$ must be an annulus
parallel to a component of $\partial X$. As the boundary of $Y$ consists of
circles which project to $\partial t$, it follows that $Y$ must project to
$E$. In this case, we homotop the two parallel components of $\partial t$
across $E$. Note that this reduces the number of components of $f^{-1}%
(f(\partial t))\cap t$ by $2$. Repeating this as needed yields the required
homotopy of $f\mid t$ to an embedding such that during the homotopy each
component of $\partial t$ either stays fixed or is moved across an annulus
component of $\partial_{1}N$.

After doing this for all components of $\partial_{0}M$, if $f\mid\partial M$
fails to be a homeomorphism, we must have two components $t_{1}$ and $t_{2}$
of $\partial_{0}M$ whose images intersect or even coincide, or we must have
components $t_{0}$ of $\partial_{0}M$ and $t_{1}$ of $\partial_{1}M$ whose
images coincide. In the first case, let $S$ denote a component of the
intersection of $t_{1}$ and $t_{2}$. (Note that $t_{1}$ and $t_{2}$ may both
be closed, in which case so is $S$.) Then $S$ is bounded by some components of
$\partial_{0}N\cap\partial_{1}N$, so that $S$ must be a union of certain
components $S_{i}$ of $\partial_{0}N$ and $\partial_{1}N$. As in the previous
paragraph, any essential loop $C$ in $S_{i}$ yields an annulus in
$(M,\partial_{0}M)$, joining $t_{1}$ and $t_{2}$, which then implies that
$\partial S_{i}$ is not empty and that $C$ is homotopic into $\partial S_{i}$.
This implies that each $S_{i}$ must be an annulus. As no component of
$\partial_{0}N$ is an annulus, it follows that $S$ must be an annulus
component of $\partial_{1}N$. Thus we can change $f$ by flip homotopies on
these annuli and obtain a homeomorphism from $\partial M$ to $\partial N$.
Finally Waldhausen's Homeomorphism Theorem shows that we can properly homotop
$f$ to a homeomorphism, completing the proof in this case. In the second case,
we have components $t_{0}$ of $\partial_{0}M$ and $t_{1}$ of $\partial_{1}M$
whose images coincide. Hence $t_{0}$ and $t_{1}$ are homotopic in $M$. If
$t_{0}$ and $t_{1}$ are closed, this implies that $M$ is homeomorphic to
$t_{0}\times I$. Otherwise, there is a submanifold $W$ of $M$ homeomorphic to
$t_{0}\times I$, such that $t_{0}\times\partial I$ corresponds to the union of
$t_{0}$ and $t_{1}$. If $C$ denotes any component of $\partial t_{0}$ and $A$
denotes the annulus in $\partial W$ joining $C$ to a component of $\partial
t_{1}$, then $A$ must be parallel into $\partial M$. It follows that $M$ is
homeomorphic to $t_{0}\times I$ in this case also. In either case, it is easy
to homotop $f$ to a homeomorphism keeping $\partial_{1}M$ mapped into
$\partial_{1}N$ during the homotopy. This completes the proof of Theorem
\ref{Jacoresult}, so long as the special case 2) in Theorem
\ref{relativeannulustheorem} never occurs.

Now suppose that for some component $t$ of $\partial_{0}M$, we are in special
case 2) of Theorem \ref{relativeannulustheorem}. As no component of
$\partial_{0}M$ is an annulus, we see that in this special case,
$(M,\partial_{0}M)$ is an $I$-bundle, whose base surface is the twice
punctured disc or the once punctured Moebius band. As $\pi_{1}(M)$ is free of
rank $2$, so is $\pi_{1}(N)$, so that each of $M$ and $N$ is a handlebody of
genus $2$.

\begin{case}
$(M,\partial_{0}M)$ is $(\Sigma\times I,\Sigma\times\partial I)$ where
$\Sigma$ denotes the twice punctured disc.
\end{case}

Thus $\partial_{1}M$ consists of three annuli. As $f$ induces a homeomorphism
$\partial_{1}M\rightarrow\partial_{1}N$, it follows that $\partial_{1}N$ also
consists of three annuli. Choosing one boundary component from each of these
three annuli in $\partial N$ gives us three disjoint simple closed curves in
$\partial N$. Each must be essential in $N$, and no two can be homotopic in
$N$, for that would yield an annulus in $M$ joining boundary circles of
distinct components of $\partial_{1}M$. As $3$ is the maximum number of
disjoint essential simple closed curves one can have on the closed orientable
surface of genus $2$, it follows that $\partial_{0}N$ must consist of two
pairs of pants, or equivalently, twice punctured discs. Now let $S$ denote one
of the two components of $\partial_{0}N$. The given map $f:(M,\partial
_{1}M)\rightarrow(N,\partial_{1}N)$ has a homotopy inverse $g:(N,\partial
_{1}N)\rightarrow(M,\partial_{1}M)$, and we consider the composite map
$\varphi:S\subset N\overset{g}{\rightarrow}M\rightarrow\Sigma$, where the last
map is simply projection. Then $\varphi:S\rightarrow\Sigma$ is proper and
$\pi_{1}$-injective and so is properly homotopic to a covering map, as $S$ is
not an annulus. As $S$ and $\Sigma$ each have Euler number equal to $-1$, it
follows that $\varphi:S\rightarrow\Sigma$ is homotopic to a homeomorphism rel
$\partial S$. Thus we can homotop $f$ to a homeomorphism while keeping
$\partial_{1}M$ mapped into $\partial_{1}N$ during the homotopy, as required.

\begin{case}
$(M,\partial_{0}M)$ is $(\Sigma\widetilde{\times}I,\Sigma\widetilde{\times
}\partial I)$ where $\Sigma$ denotes the once punctured Moebius band.
\end{case}

In this case, $\partial_{0}M$ is a connected double cover of $\Sigma$, and
$\partial_{1}M$ consists of two annuli. As before, we let $S$ denote a
component of $\partial_{0}N$. Recall that $S\;$is not an annulus. Now let
$g:(N,\partial_{1}N)\rightarrow(M,\partial_{1}M)$ denote a homotopy inverse to
$f$, and consider the composite map $\varphi:S\subset N\overset{g}%
{\rightarrow}M\rightarrow\Sigma$, where the last map is simply projection.
Then $\varphi:S\rightarrow\Sigma$ is proper and $\pi_{1}$-injective and so is
properly homotopic to a covering map, as $S$ is not an annulus. The degree of
this covering must be even, as $S$ is orientable and $\Sigma$ is not. As $S$
is an incompressible subsurface of $\partial N$, we have $\chi(S)\geq
\chi(\partial N)=-2$. As $\chi(\Sigma)=-1$, it follows that the covering has
degree $2$. It follows that each of $S$ and $\partial_{0}M$ is the orientable
double cover of $\Sigma$. Hence we can homotop $f$ to a homeomorphism while
keeping $\partial_{1}M$ mapped into $\partial_{1}N$ during the homotopy, as required.
\end{proof}

We complement the above result by observing what happens if we drop the
assumption that no component of $\partial_{0}M$ or $\partial_{0}N$ is an annulus.

\begin{theorem}
Let $(M,\partial_{0}M,\partial_{1}M)$, $(N,\partial_{0}N,\partial_{1}N)$ be
orientable $3$-manifolds with $M$ and $N$ Haken, and $\partial_{0}M$
incompressible. Suppose that any embedded incompressible annulus in
$(M,\partial_{0}M)$ is parallel into $\partial_{0}M$ or $\partial_{1}M$ and
that $M$ and $N$ are not solid tori. If $f:M\rightarrow N$ is a homotopy
equivalence such that the restriction $f\mid\partial_{1}M$ is a homeomorphism
onto $\partial_{1}N$, then $f$ can be deformed to a homeomorphism while
keeping $\partial_{1}M$ mapped into $\partial_{1}N$ during the homotopy.
\end{theorem}

\begin{proof}
As in the proof of Theorem \ref{Jacoresult}, we can homotop $f$ to be proper
while keeping $\partial_{1}M$ mapped into $\partial_{1}N$ during the homotopy.
Now we apply Theorem \ref{EvansSwarupTucker}. We conclude that $f$ is properly
homotopic to a map $f_{1}$ which is a homeomorphism unless a) $M$ is a
$I$-bundle over a closed surface and $f_{1}(M)\subset\partial N$, or b) $N$
(and hence also $M$) is a solid torus and $f_{1}$ is a branched covering with
branch set a circle, or c) $M$ is a handlebody and $f_{1}(M)\subset\partial
N$. If $f_{1}$ is a homeomorphism, we can alter $f_{1}$ and the homotopy so as
to arrange that $\partial_{1}M$ is mapped into $\partial_{1}N$ during the
homotopy, which completes the proof of the theorem in this case. Case b) is
excluded in the hypotheses of our theorem, so it remains to eliminate the two
cases a) and c).

Suppose that we have case a), so that $M$ is a $I$-bundle over a closed
surface $F$ and $f_{1}(M)$ lies in a component $S$ of $\partial N$. As $f$ is
a homotopy equivalence, the inclusion of $S$ in $N$ must also be a homotopy
equivalence, so that $N$ is homeomorphic to $S\times I$. If $M$ is the trivial
$I$-bundle over $F$, then $f$ must map each component of $\partial M$ to $S$
with degree $1$. Otherwise, the bundle is non-trivial and $f$ maps $\partial
M$ to $S$ with degree $2$. In either case, $M$ contains embedded annuli which
are not properly homotopic into $\partial M$, so $\partial_{1}M$ must be
non-empty. Let $Y$ be a component of $\partial_{1}N$ in $S$ and $X$ be a
component of $\partial_{1}M$ which maps to $Y$ by a homeomorphism. Let
$\Sigma$ denote a component of $\partial Y$, and let $C$ denote the component
of $\partial X$ mapped to $\Sigma$. Whether $M$ is a trivial or non-trivial
$I$-bundle, there is an annulus $A$ embedded in $M$ with one end equal to $C$
which cannot be properly homotoped into $\partial M$. Let $D$ denote the other
end. Suppose that $D$ meets a component $Z$ of $\partial_{1}M$. As $f$ maps
$\partial_{1}M$ to $\partial_{1}N$ by a homeomorphism, $f(Z)$ is disjoint from
$Y$. As $f(D)$ is homotopic to $f(C)=\Sigma$, we see that $f(D)$ can be
homotoped out of $f(Z)$, and hence that $D$ can be homotoped out of $Z$. By
repeating this argument for any other components of $\partial_{1}M$ which meet
$D$, we can arrange that $D$ lies in $\partial_{0}M$. But this immediately
contradicts the condition that any embedded incompressible annulus in
$(M,\partial_{0}M)$ is parallel into $\partial_{0}M$ or $\partial_{1}M$. This
contradiction shows that case a) cannot occur.

Finally suppose that we have case c), so that $M$, and hence $N$, is a
handlebody. If $t$ denotes a component of $\partial_{0}M$, the argument in the
proof of Theorem \ref{Jacoresult} shows that either $f\mid t$ can be modified
rel $\partial t$ to be an embedding, or that the image of $t$ contains an
annulus $S$ which is a union of components of $\partial_{0}M$ and
$\partial_{1}M$. If we regard $S$ as an annulus in $(M,\partial_{0}M)$, the
assumptions of our theorem imply that $S$ is parallel into $\partial_{0}M$ or
$\partial_{1}M$. The first case is only possible if $t$ is an annulus, and the
second is only possible if $S$ is a component of $\partial_{1}M$. In the
second case, we can further homotop $t$ to an embedding by moving two
components of $\partial t$ across $S$. Thus either, we can homotop $f\mid t$
to an embedding or $t$ is an annulus and the boundary of $N$ is a torus
consisting of the image of $t$. As $N$ is a handlebody, it must be a solid
torus, contradicting our assumptions. Thus $f\mid t$ can be homotoped to an
embedding relative to its boundary, for each component $t$ of $\partial_{0}M$,
and we complete the argument as in \ref{Jacoresult}.
\end{proof}

\end{document}